\documentclass{article}

\usepackage{latexsym}
\usepackage{amsfonts}
\usepackage{amssymb}
\usepackage{cite}

\usepackage{proof}

\newcommand{\dg}{\mbox{\rm dg}}

\newcommand{\blem}{\begin{lemma}}
\newcommand{\elem}{\end{lemma}}
\newcommand{\bth}{\begin{theorem}}
\newcommand{\ethm}{\end{theorem}}
\newcommand{\benu}{\begin{enumerate}}
\newcommand{\eenu}{\end{enumerate}}
\newcommand{\bdes}{\begin{description}}
\newcommand{\edes}{\end{description}}
\newcommand{\bdf}{\begin{definition}}
\newcommand{\edf}{\end{definition}}
\newcommand{\bcor}{\begin{cor}}
\newcommand{\ecor}{\end{cor}}
\newcommand{\bprp}{\begin{proposition}}
\newcommand{\eprp}{\end{proposition}}
\newcommand{\bmlem}{\begin{mlemma}}
\newcommand{\emlem}{\end{mlemma}}
\newcommand{\bclm}{\begin{claim}}
\newcommand{\eclm}{\end{claim}}
\newcommand{\bprf}{{\bf Proof}.\hspace{2mm}}

\newcommand{\eprf}{\hspace*{\fill} $\Box$}

\newcommand{\beqn}{\begin{equation}}
\newcommand{\eeqn}{\end{equation}}
\newcommand{\beqnarr}{\begin{eqnarray}}
\newcommand{\eeqnarr}{\end{eqnarray}}
\newcommand{\beqnarrs}{\begin{eqnarray*}}
\newcommand{\eeqnarrs}{\end{eqnarray*}}

\newtheorem{theorem}{Theorem}[section]
\newtheorem{definition}[theorem]{Definition}
\newtheorem{proposition}[theorem]{Proposition}
\newtheorem{lemma}[theorem]{Lemma}
\newtheorem{cor}[theorem]{Corollary}
\newtheorem{mlemma}[theorem]{Main Lemma}
\newtheorem{claim}[theorem]{Claim}

\newcommand{\Del}{\Delta}
\newcommand{\ome}{\omega}

\newcommand{\bet}{\beta}
\newcommand{\gam}{\gamma}
\newcommand{\Gam}{\Gamma}

\newcommand{\fal}{\forall}
\newcommand{\exi}{\exists}

\newcommand{\Rarw }{\Rightarrow}

\title{Proof-theoretic strengths of weak theories for positive inductive definitions}
\author{Toshiyasu Arai
\\
Graduate School of Science,
Chiba University
\\
1-33, Yayoi-cho, Inage-ku,
Chiba, 263-8522, JAPAN
}

\date{}

\begin{document}

\maketitle

\begin{abstract}
In this paper the lightface $\Pi^{1}_{1}$-Comprehension axiom is shown to be proof-theoretically strong even over 
$\mbox{{\rm RCA}}_{0}^{*}$, and we calibrate the proof-theoretic ordinals of weak 
fragments of the theory $\mbox{ID}_{1}$ of positive inductive definitions over natural numbers.
Conjunctions of negative and positive formulas in the transfinite induction axiom of $\mbox{ID}_{1}$
are shown to be weak, and disjunctions are strong.
Thus we draw a boundary line between predicatively reducible and impredicative fragments of $\mbox{ID}_{1}$.
\end{abstract}

\section{The lightface $\Pi^{1}_{1}$-Comprehension axiom over $\mbox{{\rm RCA}}^{*}_{0}$}\label{sec:intro}
This research is motivated to answer the questions raised by 
J. Van der Meeren, M. Rathjen and A. Weiermann \cite{10,11}:
Let $|T|$ denote the proof-theoretic ordinal of a subsystem $T$ of second order arithmetic.
\\
{\bf Conjecture}{\rm (\cite{10,11})}
\begin{enumerate}
\item\label{conjecture.1}
$|\mbox{{\rm RCA}}^{*}_{0}+(\Pi^{1}_{1}(\Pi^{0}_{3})\mbox{{\rm -CA}})^{-}|=\varphi\omega 0$.
\item\label{conjecture.2}
$|\mbox{{\rm RCA}}^{*}_{0}+(\Pi^{1}_{1}\mbox{{\rm -CA}})^{-}|=\varphi\omega 0$.
\item\label{conjecture.3}
$|\mbox{{\rm RCA}}_{0}+(\Pi^{1}_{1}\mbox{{\rm -CA}})^{-}|=\vartheta(\Omega^{\omega})$.
\end{enumerate}
where $\mbox{{\rm RCA}}_{0}^{*}$ defined \cite{9} is obtained from $\mbox{RCA}_{0}$
by adding a function symbol for the exponential function $2^{x}$ together with an axiom for the function $2^{x}$, 
and restricting the induction axiom schema to bounded formulas in the expanded language.
In $\mbox{RCA}_{0}$ the induction axiom schema is available for $\Sigma^{0}_{1}$-formulas.
$(\Pi^{1}_{1}(\Pi^{0}_{3})\mbox{{\rm -CA}})^{-}$ denotes the axiom schema of lightface, i.e.,
set parameter-free $\Pi^{1}_{1}$-Comprehension Axiom with $\Pi^{0}_{3}$-matrix $\lambda$:
\begin{equation}\label{eq:Pi11CA}
\exists Y\forall n[n\in Y\leftrightarrow \forall X\,\lambda(X,n)]
\end{equation}
$(\Pi^{1}_{1}\mbox{{\rm -CA}})^{-}$ is the axiom schema of set parameter-free $\Pi^{1}_{1}$-Comprehension Axiom 
with arbitrary arithmetical formulas $\lambda$.
$\varphi$ in the ordinal $\varphi\omega 0$ 
denotes the binary Veblen function, and $\vartheta$ in $\vartheta(\Omega^{\omega})$ is a collapsing function
introduced in \cite{8}.
The ordinal $\Gamma_{0}=\vartheta(\Omega^{2})$ is known to be the limit of predicativity.

When $\Sigma^{0}_{1}$-formulas are available in the induction axiom schema, the proof-theoretic ordinal is shown to
be the small Veblen ordinal $\vartheta(\Omega^{\omega})$.

\bth\label{th:JWA}{\rm (\cite{11})}
$|\mbox{{\rm RCA}}_{0}+(\Pi^{1}_{1}(\Pi^{0}_{3})\mbox{{\rm -CA}})^{-}|=\vartheta(\Omega^{\omega})$.
\end{theorem}

According to \cite{10}
A. Weiermann showed that 
the wellfoundedness of ordinals up to each ordinal$<\vartheta(\Omega^{\omega})$
is provable in $\mbox{{\rm RCA}}_{0}+(\Pi^{1}_{1}(\Pi^{0}_{3})\mbox{{\rm -CA}})^{-}$,
and M. Rathjen showed that 
$\mbox{{\rm RCA}}_{0}+(\Pi^{1}_{1}(\Pi^{0}_{3})\mbox{{\rm -CA}})^{-}$ is reducible to
 $\Pi^{1}_{2}\mbox{-BI}_{0}$,
whose proof-theoretic ordinal, $|\Pi^{1}_{2}\mbox{-BI}_{0}|$, is $\vartheta(\Omega^{\omega})$, cf.\,\cite{8}.

In trying to settle the Conjecture affirmatively,
we have first investigated weak fragments of the theory $\mbox{ID}_{1}$ of positive inductive definitions
over natural numbers, and found a line between predicatively reducible and impredicative fragments 
of $\mbox{ID}_{1}$, cf.\,Theorem \ref{th:main} below.
One fragment is proof-theoretically strong in the sense that the fragment proves the wellfoundedness up to each 
ordinal$<\vartheta(\Omega^{\omega})$, cf.\,Lemma \ref{lem:Pi0WF}.
The proof can be transformed to one in $\mbox{{\rm RCA}}_{0}^{*}+(\Pi^{1}_{1}(\Pi^{0}_{3})\mbox{{\rm -CA}})^{-}$, and thereby
we obtain
$|\mbox{{\rm RCA}}_{0}^{*}+(\Pi^{1}_{1}(\Pi^{0}_{3})\mbox{{\rm -CA}})^{-}|\geq\vartheta(\Omega^{\omega})$.
By combining Theorem \ref{th:JWA} we arrive at a negative answer to the Conjecture \ref{conjecture.1},
$|\mbox{{\rm RCA}}_{0}^{*}+(\Pi^{1}_{1}(\Pi^{0}_{3})\mbox{{\rm -CA}})^{-}|=\vartheta(\Omega^{\omega})$.
Actually in Proposition \ref{prp:mainRCA}.\ref{prp:mainRCA.1} we will see that
$\mbox{{\rm RCA}}_{0}^{*}$ is equal to $\mbox{{\rm RCA}}_{0}$ with the help of
the lightface $\Pi^{1}_{1}$-Comprehension Axiom.
Thus the whole of the Conjecture is refuted.

The theory $\mbox{ID}_{1}$ for non-iterated positive inductive definitions over natural numbers 
is an extension of the first-order arithmetic {\sf PA} in a language $\mathcal{L}(\mbox{ID})$, which is
obtained from an arithmetic language by adding unary predicate constant $R_{\varphi}$ for each 
$X$-positive formula $\varphi(X,x)$.
Axioms are
\begin{equation}\label{eq:delta0ind}
\theta(0) \land \forall x(\theta(x)\to\theta(x+1))\to \forall y\,\theta(y)
\end{equation}
for each $\mathcal{L}(\mbox{ID})$-formula $\theta$.

\begin{equation}\label{eq:Rax1}
\forall x[\varphi(R_{\varphi},x) \to R_{\varphi}(x)]
\end{equation}

\begin{equation}\label{eq:Rax2}
\forall u[R_{\varphi}(u) \to \forall x(\varphi(\sigma,x) \to \sigma(x)) \to \sigma(u)]
\end{equation}
for each $\mathcal{L}(\mbox{ID})$-formula $\sigma$.

Note that $\mbox{{\rm ID}}_{1}$
 proves that $R_{\varphi}$ is a fixed point of positive $\varphi$,
\begin{equation}\label{eq:Rax3}
\forall x[R_{\varphi}(x)\to \varphi(R_{\varphi},x)]
\end{equation}
since $\varphi(\varphi(R_{\varphi}))\subset \varphi(R_{\varphi})$ by (\ref{eq:Rax1}),
and then apply (\ref{eq:Rax2}) to the formula $\sigma(x)\equiv\varphi(R_{\varphi},x)$.

$\Pi^{0}_{k}(\Omega)\mbox{-ID}(\Pi^{0}_{1})$ is a fragment of $\mbox{ID}_{1}$ in which 
$X$-positive formulas $\varphi(X,x)\in\Pi^{0}_{1}$,
the formulas $\theta$ in the the complete induction schema (\ref{eq:delta0ind}) 
as well as
the formulas $\sigma$ in the axiom (\ref{eq:Rax2}) are restricted to $\Pi^{0}_{k}$-formulas
$\theta,\sigma\in\Pi^{0}_{k}(\Omega)$ in the language $\mathcal{L}(\mbox{ID})$ with atomic formulas $R_{\varphi}(t)$.

Let ${\sf EA}^{2}$ be the elementary recursive arithmetic in the second order logic, i.e.,
no comprehension axiom such as $\Delta^{0}_{0}$-Comprehension is assumed.
IND denotes the $\Pi^{1}_{1}$-sentence
$\forall X\forall a[X(0)\land \forall y(X(y)\to X(y+1))\to X(a)]$.
For a first-order positive formula $\varphi(X,y)$, let
\[
I_{\varphi}=\bigcap\{X: \varphi(X)\subset X\}.
\]
The following Proposition \ref{prp:mainRCA} is utilized to show Theorem \ref{th:mainRCA}, in which the base theory ${\sf EA}^{2}$
can be the stronger $\mbox{{\rm RCA}}_{0}^{*}$.

\begin{proposition}\label{prp:mainRCA}
Let $k\geq 1$.
\begin{enumerate}
\item\label{prp:mainRCA.1}
${\sf EA}^{2}+\mbox{{\rm IND}}+(\Pi^{1}_{1}(\Sigma^{0}_{k+1})\mbox{{\rm -CA}})^{-}\vdash\Sigma^{0}_{k}\mbox{{\rm -IND}}$.
\item\label{prp:mainRCA.2}
Let $\varphi(X,y)\in\Pi^{0}_{1}$ be an $X$-positive formula.
If $\Pi^{0}_{k}(\Omega)\mbox{{\rm -ID}}(\Pi^{0}_{1})\vdash A(R_{\varphi})$, then
${\sf EA}^{2}+\mbox{{\rm IND}}+(\Pi^{1}_{1}(\Sigma^{0}_{k+1})\mbox{{\rm -CA}})^{-}\vdash A(I_{\varphi})$.
\item\label{prp:mainRCA.3}
For any $X$-positive formula $\varphi(X,y)$,
if $\mbox{{\rm ID}}_{1}\vdash A(R_{\varphi})$, then
${\sf EA}^{2}+\mbox{{\rm IND}}+(\Pi^{1}_{1}\mbox{{\rm -CA}})^{-}\vdash A(I_{\varphi})$.
\end{enumerate}
\end{proposition}
\bprf
Let $k\geq 1$.
\\
\ref{prp:mainRCA}.\ref{prp:mainRCA.1}.
For a $\Sigma^{0}_{k}$-formula $\varphi(a,X,z)$ let
\[
N(a,z) :\Leftrightarrow \forall X[\varphi(0,X,z) \land\forall y(\varphi(y,X,z) \to\varphi(y+1,X,z))\to \varphi(a,X,z)]
\]
$N(a,z)$ is a $\Pi^{1}_{1}(\Sigma^{0}_{k+1})$-formula without set parameter, 
and exists as a set by the axiom $(\Pi^{1}_{1}(\Sigma^{0}_{k+1})\mbox{{\rm -CA}})^{-}$.
It is obvious that $N(a,z)$ is inductive with respect to $a$, i.e., $N(0,z)$ and $\forall a[N(a,z)\to N(a+1,z)]$.
Therefore by IND we obtain $\forall a\, N(a,z)$, i.e.,
$\forall X[\varphi(0,X,z) \land\forall y(\varphi(y,X,z) \to\varphi(y+1,X,z))\to \forall a\,\varphi(a,X,z)]$.
\\

\noindent
\ref{prp:mainRCA}.\ref{prp:mainRCA.2} and \ref{prp:mainRCA}.\ref{prp:mainRCA.3}.
We show Proposition \ref{prp:mainRCA}.\ref{prp:mainRCA.2}. Proposition \ref{prp:mainRCA}.\ref{prp:mainRCA.3} is similarly seen.
Argue in ${\sf EA}^{2}+\mbox{{\rm IND}}+(\Pi^{1}_{1}(\Sigma^{0}_{k+1})\mbox{{\rm -CA}})^{-}$, and 
let $\varphi,\psi$ be
positive $\Pi^{0}_{1}$-formulas.
First note that $I_{\varphi}$ exists as a set by $(\Pi^{1}_{1}(\Sigma^{0}_{2})\mbox{{\rm -CA}})^{-}$.
By Proposition \ref{prp:mainRCA}.\ref{prp:mainRCA.1} we have
$\Sigma^{0}_{k}\mbox{{\rm -IND}}$, and hence
$A(I_{\varphi},0)\land\forall x(A(I_{\varphi},x)\to A(I_{\varphi},x+1))\to \forall z\,A(I_{\varphi},z)$ for any 
$\Sigma^{0}_{k}$- or $\Pi^{0}_{k}$-formula $A$.
$\varphi(I_{\varphi})\subset I_{\varphi}$ is seen logically.

Let $A(I_{\psi},y)$ be a $\Pi^{0}_{k}$-formula. We show
\begin{equation}\label{eq:AIvphi}
\varphi(A(I_{\psi}))\subset A(I_{\psi}) \to I_{\varphi}\subset A(I_{\psi})
\end{equation}
Let 
\[
B(y):\Leftrightarrow \forall X[\varphi(A(X))\subset A(X) \to A(X,y)]
.\]
$B$ exists as a set by $(\Pi^{1}_{1}(\Sigma^{0}_{k+1})\mbox{{\rm -CA}})^{-}$.
We claim that 
\begin{equation}\label{eq:Bcl}
\varphi(B)\subset B
\end{equation}
Assume $\varphi(B,y)$ and $\varphi(A(X))\subset A(X)$.
We need to show $A(X,y)$.
We first show $B\subset A(X)$.
Suppose $B(z)$. Then by the assumption $\varphi(A(X))\subset A(X)$ we have $A(X,z)$.
Hence $B\subset A(X)$, and $\varphi(B,y)\to \varphi(A(X),y)$ by the positivity of $\varphi(X)$.
The assumption $\varphi(B,y)$ yields $\varphi(A(X),y)$, and we conclude $A(X,y)$ by $\varphi(A(X))\subset A(X)$.

Since $B$ is a set, we obtain $I_{\varphi}\subset B$ by (\ref{eq:Bcl}).
On the other hand we have $B(y)\to \varphi(A(I_{\psi}))\subset A(I_{\psi}) \to A(I_{\psi},y)$ since $I_{\psi}$
is a set.
Therefore $\varphi(A(I_{\psi}))\subset A(I_{\psi})\to B\subset A(I_{\psi})$.
This together with $I_{\varphi}\subset B$ yields (\ref{eq:AIvphi}).
\eprf

Let
$\omega_{0}(\alpha)=\alpha$, and $\omega_{n+1}(\alpha)=\omega^{\omega_{n}(\alpha)}$ for $n\geq 0$.
Thus $\Omega^{\omega}=\omega_{1}(\Omega\cdot \omega)$.

\bth\label{th:mainRCA}
Let $k\geq 1$.
\begin{enumerate}
\item\label{th:mainRCA.1}
$\mbox{{\rm RCA}}_{0}^{*}+(\Pi^{1}_{1}(\Pi^{0}_{k+2})\mbox{{\rm -CA}})^{-}=
\mbox{{\rm RCA}}_{0}+(\Pi^{1}_{1}(\Pi^{0}_{k+2})\mbox{{\rm -CA}})^{-}$.
\item\label{th:mainRCA.2}
$\Pi^{0}_{k}(\Omega)\mbox{{\rm -ID}}(\Pi^{0}_{1})$ is interpreted canonically in
$\mbox{{\rm RCA}}_{0}+(\Pi^{1}_{1}(\Pi^{0}_{k+2})\mbox{{\rm -CA}})^{-}$.
\begin{eqnarray*}
&& |\Pi^{0}_{k}(\Omega)\mbox{{\rm -ID}}(\Pi^{0}_{1})| = |\mbox{{\rm RCA}}^{*}_{0}+(\Pi^{1}_{1}(\Pi^{0}_{k+2})\mbox{{\rm -CA}})^{-}|
\\
& = &
|\mbox{{\rm RCA}}_{0}+(\Pi^{1}_{1}(\Pi^{0}_{k+2})\mbox{{\rm -CA}})^{-}|
= |\Pi^{1}_{k+1}\mbox{{\rm -BI}}_{0}|=\vartheta(\omega_{k}(\Omega\cdot\omega))
\end{eqnarray*}

\item\label{th:mainRCA.3}
$\mbox{{\rm ID}}_{1}$ is interpreted canonically in
$\mbox{{\rm RCA}}_{0}+(\Pi^{1}_{1}\mbox{{\rm -CA}})^{-}$.
\[
|\mbox{{\rm ID}}_{1}|=|\mbox{{\rm RCA}}^{*}_{0}+(\Pi^{1}_{1}\mbox{{\rm -CA}})^{-}|=|\mbox{{\rm RCA}}_{0}+(\Pi^{1}_{1}\mbox{{\rm -CA}})^{-}|=\vartheta(\varepsilon_{\Omega+1})
.\]
\end{enumerate}
\end{theorem}
\bprf
Theorem \ref{th:mainRCA}.\ref{th:mainRCA.1} follows from Proposition \ref{prp:mainRCA}.\ref{prp:mainRCA.1}.

Let us consider Theorem \ref{th:mainRCA}.\ref{th:mainRCA.2}.
From  Proposition \ref{prp:mainRCA}.\ref{prp:mainRCA.2} we see that 
$\Pi^{0}_{k}(\Omega)\mbox{{\rm -ID}}(\Pi^{0}_{1})$ is interpreted canonically in
$\mbox{{\rm RCA}}_{0}+(\Pi^{1}_{1}(\Pi^{0}_{k+2})\mbox{{\rm -CA}})^{-}$.
As in \cite{11}
using \cite{6}
we see that $\Pi^{1}_{k+1}\mbox{{\rm -BI}}_{0}$ comprises
$\mbox{{\rm RCA}}_{0}+(\Pi^{1}_{1}(\Pi^{0}_{k+2})\mbox{{\rm -CA}})^{-}$.
In \cite{8}
it is shown that $|\Pi^{1}_{k+1}\mbox{{\rm -BI}}_{0}|=\vartheta(\omega_{k}(\Omega\cdot\omega))$.
The fact $\vartheta(\omega_{k}(\Omega\cdot\omega))\leq|\Pi^{0}_{k}(\Omega)\mbox{{\rm -ID}}(\Pi^{0}_{1})|$ is shown in
Proposition \ref{prp:lowerb} below.
\eprf

Therefore even the weakest fragment $\mbox{{\rm RCA}}_{0}^{*}+(\Pi^{1}_{1}(\Pi^{0}_{3})\mbox{{\rm -CA}})^{-}$
is not predicatively reducible.
In \cite{10}
it is reported that 
$|\mbox{{\rm RCA}}_{0}^{*}+(\Pi^{1}_{1}(\Pi^{0}_{2})\mbox{{\rm -CA}})^{-}|=|I\Sigma_{2}|$
and $|\mbox{{\rm RCA}}_{0}+(\Pi^{1}_{1}(\Pi^{0}_{2})\mbox{{\rm -CA}})^{-}|=|I\Sigma_{3}|$.

This indicates that fragments of the light-face $\Pi^{1}_{1}$-CA could be graded
according to another hierarchy of formulas for matrices $\lambda$ in (\ref{eq:Pi11CA})
rather than the usual arithmetic hierarchy.
Let $\Pi^{0}_{0}(\mbox{P}^{X})=\Sigma^{0}_{0}(\mbox{P}^{X})$ denote the class of first-order formulas 
$\lambda(X)$ obtained from arithmetic atomic formulas and $X$-positive formulas
by means of Boolean connectives and bounded quantifications.
Classes $\Pi^{0}_{k}(\mbox{P}^{X})$ and $\Sigma^{0}_{k}(\mbox{P}^{X})$ of first-order formulas are defined
from the class by prefixing alternating (unbounded) quantifiers.
It is open for us, but seems to me plausible that
$\mbox{{\rm RCA}}^{*}_{0}+(\Pi^{0}_{2}(\mbox{{\rm P}}^{X})\mbox{{\rm -CA}})^{-}$ is, or even
$\mbox{{\rm RCA}}_{0}+(\Pi^{0}_{2}(\mbox{{\rm P}}^{X})\mbox{{\rm -CA}})^{-}$
is predicatively reducible.

\subsection{Weak fragments}
Let us introduce weak fragments of the theory $\mbox{ID}_{1}$.

Let $\mathcal{L}$ be a language for arithmetic having function constants\footnote{The proof-theoretic strength does not increase with more constants, e.g., with function constants for primitive recursive functions.}
for each elementary recursive functions.
Relation symbols in $\mathcal{L}$ are $=,<$.
$\Delta_{0}^{-}$ denotes the set of bounded formulas in $\mathcal{L}$, and $\Pi^{1-}_{0}$
the set of formulas in $\mathcal{L}$, called \textit{arithmetical} formulas.
The elementary recursive arithmetic {\sf EA} is the theory in $\mathcal{L}$ whose axioms are
defining axioms for function constants, axioms for $=,<$ and $\Delta_{0}^{-}\mbox{-IND}$:
(\ref{eq:delta0ind}) is restricted to $\theta\in\Delta_{0}^{-}$.

For a second-order arithmetic $T$, its \textit{proof-theoretic ordinal} $|T|$ is defined to be
the supremum of the order types $|\prec|$ of elementary recursive and transitive relations $\prec$ for which
$T\vdash \forall y(\forall x\prec y\, X(x)\to X(y))\to\forall y\,X(y)$.
When $T$ is a theory for positive inductive definitions,
$|T|$ is defined to be 
the supremum of the order types $|\prec|$ of elementary recursive and transitive relations $\prec$ for which
$T\vdash \forall x(x\in W_{\prec})$ for the accessible (well founded) part $W_{\prec}$ of the relation $\prec$.

For a class $\Phi$ of $X$-positive formulas $\varphi(X,x)$,
let $\mathcal{L}(\Phi)=\mathcal{L}\cup\{R_{\varphi}: \varphi\in\Phi\}$ denote the language obtained from the language $\mathcal{L}$
by adding unary predicate constants $R_{\varphi}$ for each $\varphi\in\Phi$.
The unary predicate constant $R_{\varphi}$ is intended to denote the least fixed point of the monotone operator defined from $\varphi$:
$\mathbb{N}\supset \mathcal{X}\mapsto\{ n\in\mathbb{N} : \mathbb{N}\models\varphi[\mathcal{X},n]\}$.

For classes $\Theta,\Gamma$ of formulas in $\mathcal{L}(\Phi)$,
let $(\Theta,\Gamma)\mbox{-ID}(\Phi)$ denote the fragment of $\mbox{ID}_{1}$ defined as follows.
$(\Theta,\Gamma)\mbox{-ID}(\Phi)$ is an extension of {\sf EA}.
In $(\Theta,\Gamma)\mbox{-ID}(\Phi)$, the positive formula $\varphi$ is in $\Phi$,
the formulas $\theta$ in
complete induction schema (\ref{eq:delta0ind}) 
are in $\Theta$,
and the formulas $\sigma$ in the axiom (\ref{eq:Rax2}) are in $\Gamma$.

When $\Theta=\Gamma$, let us write $\Gamma\mbox{-ID}(\Phi)$ for $(\Gamma,\Gamma)\mbox{-ID}(\Phi)$, and
when $\Phi$ is the class 
of all positive formulas, let us write $\Gamma\mbox{-ID}$ for $\Gamma\mbox{-ID}(\Phi)$.
\\

${\sf Acc}$ denotes the class of $X$-positive formulas $\varphi$
 \begin{equation}\label{eq:Acc}
 \varphi(X,x) \equiv  [\forall y(\theta_{0}(x,y)\to t_{0}(x,y)\in X)]
 \end{equation}
with an arithmetic bounded formula $\theta_{0}(x,y)$ and a term $t_{0}(x,y)$.
In $\theta_{0}(x,y)$ and $t_{0}(x,y)$ first-order parameters other than $x,y$ may occur.
 For an elementary recursive relation $\prec$,
  $\forall y(y\prec x\to y\in X)$ is a typical example of an ${\sf Acc}$-operator.
$Acc$ denotes the class of formulas $\sigma(x)$ which are obtained from an {\sf Acc}-operator 
 by substituting any predicate constant $R$ for $X$
\begin{equation}\label{eq:sigma}
 \sigma(x)\equiv [\forall y(\theta_{0}(x,y)\to t_{0}(x,y)\in R)]
\end{equation}
where $\theta_{0}(x,y)$ is  an arithmetic bounded formula and $t_{0}(x,y)$ a term  possibly with 
first-order parameters other than $x,y$.

\bdf\label{df:PN}
{\rm
A formula 
 is said to be \textit{positive} [\textit{negative}] if each predicate constant $R_{\varphi}$ for least fixed point
occurs only positively [negatively] in it, resp.
$\mbox{Pos}$ [$\mbox{Neg}$] denotes the class of all positive formulas [the class of all negative formulas], resp.

Also let $\mbox{P}\cup\mbox{N}:=\mbox{Pos}\cup\mbox{Neg}$,
$\mbox{P}\land\mbox{N}:=\{C\land D: C\in\mbox{Pos}, D\in\mbox{Neg}\}$ and
$\mbox{N}\lor\mbox{P}:=\{D\lor C : D\in\mbox{Neg}, C\in\mbox{Pos}\}$.}
\edf

\bdf\label{df:SigPi}
{\rm For $k\geq 0$, classes $\Pi^{0}_{k}(\mbox{{\rm P}})$ and $\Sigma^{0}_{k}(\mbox{{\rm P}})$ of formulas in the language 
$\mathcal{L}(\mbox{{\rm ID}})$ are defined recursively.
\begin{enumerate}
\item
$\Pi^{0}_{0}(\mbox{{\rm P}})=\Sigma^{0}_{0}(\mbox{{\rm P}})$ denotes the class of bounded formulas in positive formulas.
Each formula in $\Pi^{0}_{0}(\mbox{{\rm P}})$ is obtained from positive formulas by means of propositional connectives
$\lnot,\lor,\land$ and bounded quantifiers $\exists x<t,\forall x<t$.

\item
$\Pi^{0}_{k}(\mbox{{\rm P}})\cup\Sigma^{0}_{k}(\mbox{{\rm P}})\subset\Pi^{0}_{k+1}(\mbox{{\rm P}})\cap \Sigma^{0}_{k+1}(\mbox{{\rm P}})$.

\item
Each class $\Pi^{0}_{k}(\mbox{{\rm P}})$ and $\Sigma^{0}_{k}(\mbox{{\rm P}})$  is closed under positive boolean combinations 
$\lor,\land$ and bounded quantifications.
\item
If $A\in \Pi^{0}_{k}(\mbox{{\rm P}})$ [$A\in \Sigma^{0}_{k}(\mbox{{\rm P}})$], then $\lnot A\in \Sigma^{0}_{k}(\mbox{{\rm P}})$
[$\lnot A\in \Pi^{0}_{k}(\mbox{{\rm P}})$], resp.

\item
If $A\in \Pi^{0}_{k}(\mbox{{\rm P}})$ [$A\in \Sigma^{0}_{k}(\mbox{{\rm P}})$], then $\forall x A\in \Pi^{0}_{k}(\mbox{{\rm P}})$
[$\exists x A\in \Sigma^{0}_{k}(\mbox{{\rm P}})$], resp.
\end{enumerate}
Let $ \Pi^{0}_{\infty}(\mbox{{\rm P}})=\bigcup_{k<\omega} \Pi^{0}_{k}(\mbox{{\rm P}})$.

Classes $\Pi^{0}_{k}(\Omega)$ and $\Sigma^{0}_{k}(\Omega)$ are defined similarly by letting
$\Pi^{0}_{0}(\Omega)=\Sigma^{0}_{0}(\Omega)$ denote the class of bounded formulas in atomic formulas $t\in R_{\varphi}, t=s,t<s$.
The predicates $R_{\varphi}$ may occur positively and/or negatively in $\Pi^{0}_{0}(\mbox{{\rm P}})$-formulas
and in $\Pi^{0}_{0}(\Omega)$-formulas.

}
\edf

$\mbox{P-ID}$ denotes the theory, in which
$X$-positive formulas $\varphi(X,x)$ are arbitrary,
the formulas $\theta$ in the complete induction schema (\ref{eq:delta0ind}) 
as well as
the formulas $\sigma$ in the axiom (\ref{eq:Rax2}) are restricted to positive formulas
$\theta,\sigma\in \mbox{Pos}$.
$(\Pi^{0}_{\infty}(\mbox{{\rm P}}),\mbox{P})\mbox{-ID}$ is an extension of $\mbox{P-ID}$ in which 
the formulas $\theta$ in (\ref{eq:delta0ind}) are arbitrary, but
$\sigma\in \mbox{Pos}$.

The following theorem is shown by D. Probst\cite{7},
and independently by B. Afshari and M. Rathjen\cite{1}.
$\mbox{P-ID}$ [$(\Pi^{0}_{\infty}(\mbox{{\rm P}}),\mbox{P})\mbox{-ID}$] is denoted as $\mbox{ID}_{1}^{*}\!\upharpoonright\!$ [$\mbox{ID}_{1}^{*}$] in \cite{7},
and as
$\mbox{ID}_{1}^{*}$ [$\mbox{ID}_{1}^{*}+\mbox{IND}_{\mathbb{N}}$] in \cite{1},
resp.

\bth\label{th:PID}{\rm (Probst\cite{7}, Afshari and Rathjen\cite{1})}
\begin{enumerate}
\item
$|\mbox{{\rm P-ID}}|=\varphi\omega 0=\vartheta(\Omega\cdot\omega)$.
\item
$|(\Pi^{0}_{\infty}(\mbox{{\rm P}}),\mbox{{\rm P}})\mbox{{\rm -ID}}|=\varphi\varepsilon_{0}0=\vartheta(\Omega\cdot\varepsilon_{0})$.
\end{enumerate}
\end{theorem}
They show that
 $\mbox{{\rm P-ID}}$ [$(\Pi^{0}_{\infty}(\mbox{{\rm P}}),\mbox{{\rm P}})\mbox{{\rm -ID}}$]
is interpreted in $\Sigma^{1}_{1}\mbox{-DC}_{0}$ [in $\Sigma^{1}_{1}\mbox{-DC}$], resp.
It seems that their proofs do not work when $\sigma$ in the axiom (\ref{eq:Rax2})
is a negative formula.

On the other side G. J\"ager and T. Strahm\cite{5} show directly the following.
Let $\mbox{ID}_{1}^{\#}$ be a subtheory of
$\widehat{ID}_{1}$ for fixed points with the axioms (\ref{eq:Rax1}) and (\ref{eq:Rax3}).
In $\mbox{ID}_{1}^{\#}$ the complete induction schema (\ref{eq:delta0ind}) is restricted to
positive formulas $\theta\in\mbox{POS}$.
\bth\label{th:JS96}{\rm (J\"ager and Strahm\cite{5})}
$|\mbox{ID}_{1}^{\#}|=\varphi\omega 0$.
\end{theorem}

In this paper we show the following theorem \ref{th:main}.
$(Acc,\mbox{{\rm N}}\lor\mbox{{\rm P}})\mbox{{\rm -ID}}({\sf Acc})$
 is the theory, in which
$X$-positive formulas $\varphi(X,x)$ are restricted to ${\sf Acc}$-operators (\ref{eq:Acc}), 
the formulas $\theta$ in the complete induction schema (\ref{eq:delta0ind}) 
are restricted to 
$\theta\in Acc$ (\ref{eq:sigma}), and
the formulas $\sigma$ in the axiom (\ref{eq:Rax2}) are restricted to a disjunction of negative formula and a positive formula
$\sigma\in \mbox{N}\lor\mbox{P}$.

$(\Pi^{0}_{k}(\mbox{{\rm P}}),\mbox{{\rm P}}\cup\mbox{{\rm N}})\mbox{{\rm -ID}}$
denotes the theory, in which
$\varphi(X,x)$ is an arbitrary $X$-positive formula,
$\theta\in \Pi^{0}_{k}(\mbox{{\rm P}})$ in (\ref{eq:delta0ind}),
and $\sigma\in \mbox{P}\cup\mbox{N}$ in (\ref{eq:Rax2}).

$(\Pi^{0}_{k}(\mbox{{\rm P}}),\mbox{{\rm P}}\land\mbox{{\rm N}})\mbox{{\rm -ID}}({\sf Acc})$ is the theory, in which
$\varphi(X,x)$ are restricted to ${\sf Acc}$-operators (\ref{eq:Acc}), 
$\theta\in \Pi^{0}_{k}(\mbox{{\rm P}})$ in (\ref{eq:delta0ind}), and
 $\sigma$ in (\ref{eq:Rax2}) are restricted to a conjunction of positive formula and a negative formula
$\sigma\in \mbox{P}\land\mbox{N}$.
$\mbox{P-ID}\subset (\Pi^{0}_{0}(\mbox{{\rm P}}),\mbox{{\rm P}}\cup\mbox{{\rm N}})\mbox{{\rm -ID}}$
is obvious.

Let $\omega_{0}=1$ and $\omega_{n+1}=\omega^{\omega_{n}}$.

\bth\label{th:main} 
\begin{enumerate}
\item\label{th:maina}
\begin{eqnarray*}
|(Acc,\mbox{{\rm N}}\lor\mbox{{\rm P}})\mbox{{\rm -ID}}({\sf Acc})|
& = & |(Acc,\Pi^{0}_{0}(\mbox{{\rm P}}))\mbox{{\rm -ID}}({\sf Acc})|
= |\Pi^{0}_{1}(\mbox{{\rm P}})\mbox{{\rm -ID}}({\sf Acc})|
\\
& = & \vartheta(\Omega^{\omega})
\end{eqnarray*}

\item\label{th:mainb}
\begin{eqnarray*}
|Acc\mbox{{\rm -ID}}({\sf Acc})| & = & |\mbox{{\rm P-ID}}|=
|(\Pi^{0}_{0}(\mbox{{\rm P}}),\mbox{{\rm P}}\cup\mbox{{\rm N}})\mbox{{\rm -ID}}|
\\
& = & |(\Pi^{0}_{0}(\mbox{{\rm P}}),\mbox{{\rm P}}\land\mbox{{\rm N}})\mbox{{\rm -ID}}({\sf Acc})|=
\vartheta(\Omega\cdot\omega)=\varphi\omega 0
\end{eqnarray*}

\item\label{th:mainc}
For each $k> 0$
\begin{eqnarray*}
&& |(\Pi^{0}_{k}(\mbox{{\rm P}}),Acc)\mbox{{\rm -ID}}({\sf Acc})| =
|(\Pi^{0}_{k}(\mbox{{\rm P}}),\mbox{{\rm P}}\cup\mbox{{\rm N}})\mbox{{\rm -ID}}|
\\
& = &
 |(\Pi^{0}_{k}(\mbox{{\rm P}}),\mbox{{\rm P}}\land\mbox{{\rm N}})\mbox{{\rm -ID}}({\sf Acc})|
=
 \vartheta(\Omega\cdot\omega_{1+k})
=\varphi\omega_{1+k}0
\end{eqnarray*}
\end{enumerate}
In particular
$|(\Pi^{0}_{\infty}(\mbox{{\rm P}}),Acc)\mbox{{\rm -ID}}({\sf Acc})|=|(\Pi^{0}_{\infty}(\mbox{{\rm P}}),\mbox{{\rm P}})\mbox{{\rm -ID}}|=
|(\Pi^{0}_{\infty}(\mbox{{\rm P}}),\mbox{{\rm P}}\cup\mbox{{\rm N}})\mbox{{\rm -ID}}|
=|(\Pi^{0}_{\infty}(\mbox{{\rm P}}),\mbox{{\rm P}}\land\mbox{{\rm N}})\mbox{{\rm -ID}}({\sf Acc})|)=\vartheta(\Omega\cdot\varepsilon_{0})=\varphi\varepsilon_{0} 0$.
\end{theorem}
Among other things this means that negative formulas $\sigma$ in the axiom (\ref{eq:Rax2}) does not raise
the proof-theoretic ordinals.
Theorem \ref{th:main}.\ref{th:mainb} strengthens Theorems \ref{th:PID} in 
\cite{7, 1}, and Theorem \ref{th:JS96} in \cite{5}.
Our proof of the upper bound is directly done by cut-eliminations in infinitary derivations.
\\

Let us mention the contents of the paper.
In Section \ref{sect:wfproof} the easy halves in Theorem \ref{th:main} are shown by giving
some wellfoundedness proofs.
In Section \ref{sec:finitarycal} theories to be considered are reformulated in one-sided sequent calculi.
In Section \ref{sec:infderiv} finitary proofs in sequent calculi are first embedded to infinitary derivations to eliminate
cut inferences partially.
This first step is needed to unfold complex induction formulas.
Second finitary proofs and infinitary derivations are embedded into a system with the operator controlled derivations
 due to W. Buchholz\cite{3}.
In the latter derivations, cut formulas are restricted to boolean combinations of positive formulas.
The upper bounds of the proof-theoretic ordinals are obtained through collapsing and bounding lemmas.
Finally we conclude the other halves in Theorem \ref{th:main}.

\section{Wellfoundedness proofs}\label{sect:wfproof}
In this section the easy halves in Theorem \ref{th:main} are shown by giving
some wellfoundedness proofs.
First let us recall the notation system $OT^{\prime}(\vartheta)$ in \cite{10}.
$OT^{\prime}(\vartheta)$ denotes a notation system of ordinals based on symbols $\{0,\Omega,+,\vartheta\}$.
 \begin{enumerate}
 \item
 $0\in OT^{\prime}(\vartheta)$. $0$ is the least element in $OT^{\prime}(\vartheta)$, and $K(0)=\emptyset$.
 \item
 If $\{\beta_{k},\alpha_{k}: k<n\}\subset OT^{\prime}(\vartheta)$ with $n>0$, $\alpha_{n-1}>\cdots>\alpha_{0}$,
 $0\neq\beta_{k}<\Omega$
 and $\alpha_{0}>0\lor n>1$,
 then $\Omega^{\alpha_{n-1}}\beta_{n-1}+\cdots+\Omega^{\alpha_{0}}\beta_{0}\in OT^{\prime}(\vartheta)$.
 $K(\Omega^{\alpha_{n-1}}\beta_{n-1}+\cdots+\Omega^{\alpha_{0}}\beta_{0})=\bigcup\{K(\alpha_{k})\cup\{\beta_{k}\}: k<n\}$.
 \item
 If $\beta\in OT^{\prime}(\vartheta)$, then $\vartheta(\beta)\in OT^{\prime}(\vartheta)\cap\Omega$.
 $K(\vartheta(\beta))=\{\vartheta(\beta)\}$.
 \item
 $
\vartheta(\alpha)<\vartheta(\beta) \Leftrightarrow 
[\alpha<\beta\land \forall\gamma\in K(\alpha)(\gamma<\vartheta(\beta))]\lor[\exists\delta\in K(\beta)(\vartheta(\alpha)\leq \delta)]
$.
\item
Each ordinal $\vartheta(\alpha)$ is defined to be additively closed.
This means that 
$\beta,\gamma<\vartheta(\alpha)\Rightarrow \beta+\gamma<\vartheta(\alpha)$.
 \end{enumerate}
Note that the system $OT^{\prime}(\vartheta)$ is $\omega$-exponential-free
except $\vartheta(\alpha)=\omega^{\alpha_{0}}$ for some $\alpha_{0}$.
An inspection of the proof in \cite{11}
shows that $Acc\mbox{-ID}({\sf Acc})$ suffices to prove 
 the wellfoundedness of ordinals up to each ordinal$<\vartheta(\Omega\cdot\omega)$. 

Let $<$ be the elementary recursive relation obtained from the relation $<$ on $OT^{\prime}(\vartheta)$
through a suitable encoding.
For the formula  $\forall y(y< x\to y\in X)$ in {\sf Acc}, let $W$ denote the accessible part of $<$, and
$Prog(X):\Leftrightarrow \forall \alpha[\forall\beta<\alpha(\beta\in X) \to \alpha\in X]$.
Then the axiom (\ref{eq:Rax1}) states $Prog(W)$,
and the axiom (\ref{eq:Rax2}) runs
$\forall x[x\in W \to Prog(\sigma) \to \sigma(x)]$ for $\sigma\in Acc$.

The following lemma shows the easy half in Theorem \ref{th:main}.\ref{th:mainb}.

\blem\label{lem:wfprf}
$Acc\mbox{{\rm -ID}}({\sf Acc})\vdash \forall\beta<\vartheta(\Omega\cdot k)(\beta\in W)$ for {\rm each} $k<\omega$.

\elem
\bprf
We see that the following are provable in $Acc\mbox{-ID}({\sf Acc})$.
Note that $A,B,C\in Acc$ for the formulas $A,B,C$ below.

\begin{enumerate}
\item[(a)]
$x\in W\to \forall y(y< x\to y\in W)$ by the axiom (\ref{eq:Rax2})
for the formula $A(x)\Leftrightarrow  \forall y(y< x\to y\in W)$.

\item[(b)]
$y\in W \to x\in W \to x+y\in W$ by the axiom (\ref{eq:Rax2})
for the formula $B(y)\Leftrightarrow  (x+y\in W)$.

\item[(c)]
Assume $K(a)\subset W$, $\forall\beta<_{\Omega}\Omega\cdot a(\vartheta(\beta)\in W)$ and $\zeta\in W\cap\Omega$, where
$\alpha<_{\Omega}\beta:\Leftrightarrow (K(\alpha)\cup K(\beta)\subset W \land \alpha<\beta)$.
Then $Prog(C)$ for $C(\xi):\Leftrightarrow (\xi<\zeta \to \vartheta(\Omega\cdot a+\xi)\in W)$.

Suppose $\xi<\zeta$ and $\forall\eta<\xi\, C(\eta)$. Then $\xi\in W$ by $\zeta\in W$, and $K(\xi)\subset W$. 
We show $\forall\alpha<\vartheta(\Omega\cdot a+\xi)(\alpha\in W)$ by $Acc$-induction on the length of $\alpha$.
By (b) we can assume that $\alpha=\vartheta(\beta)$.
If $\alpha\leq \xi_{0}$ for a $\xi_{0}\in K(a)\cup K(\xi)\subset W$, then $\alpha\in W$.
Otherwise $K(\beta)\subset W$ by the induction hypothesis, and $\beta<\Omega\cdot a+\xi$.
We can assume $\beta=\Omega\cdot a+\eta$ for an $\eta<\xi$ by 
the assumption $\forall\beta<_{\Omega}\Omega\cdot a(\vartheta(\beta)\in W)$.
We obtain $\alpha\in W$ by $C(\eta)$.

\item[(d)]
$K(a)\subset W \to \forall\beta<_{\Omega}\Omega\cdot a(\vartheta(\beta)\in W) \to \forall\beta<_{\Omega}\Omega\cdot(a+1)(\vartheta(\beta)\in W)$.

Assume $K(a)\subset W$, $\forall\beta<_{\Omega}\Omega\cdot a(\vartheta(\beta)\in W)$ and $\beta<_{\Omega}\Omega\cdot(a+1)$.
We need to show $\vartheta(\beta)\in W$.
We can assume $\beta=\Omega\cdot a+\zeta$ for a $\zeta<\Omega$.
By $K(\beta)\subset W$ we have $\zeta\in W$.
Then by (c) we have $Prog(C)$, which yields $\forall\xi\in W\cap\zeta(\vartheta(\Omega\cdot a+\xi)\in W)$
by the axiom (\ref{eq:Rax2}) for the $Acc$-formula $C$.
From this we see that
$\forall\alpha<\vartheta(\Omega\cdot a+\zeta)(\alpha\in W)$ by $Acc$-induction on the length of $\alpha$, and hence
$\vartheta(\Omega\cdot a+\zeta)\in W$ as desired.

\end{enumerate}
By (d) we obtain $\forall\beta<_{\Omega}\Omega\cdot k\, C(\beta)$, i.e.,
$\forall\beta<_{\Omega}\Omega\cdot k(\vartheta(\beta)\in W)$ for each $k<\omega=\vartheta(1)$ with $K(k)=\{1\}=\{\vartheta(0)\}\subset W$.
Using this and (b), we see that $\forall\beta<\vartheta(\Omega\cdot k)(\beta\in W)$
by $Acc$-induction on the length of $\beta$.
This shows Lemma \ref{lem:wfprf}.
\eprf

Although the fact $\vartheta(\omega_{k}(\Omega\cdot\omega))\leq|\Pi^{0}_{k}(\Omega)\mbox{{\rm -ID}}(\Pi^{0}_{1})|$ is assumed to be a folklore,
let us give a proof of it for completeness, cf.\,Theorem \ref{th:mainRCA}.\ref{th:mainRCA.2}.
Let $\Omega_{0}(\alpha)=\alpha$, and $\Omega_{n+1}(\alpha)=\Omega^{\Omega_{n}(\alpha)}$.
Then $\omega_{k}(\Omega\cdot\omega)=\Omega_{k}(\omega)$ for $k\geq 1$.

\bprp\label{prp:lowerb}
Let $k\geq 1$.
$\Pi^{0}_{k}(\Omega)\mbox{{\rm -ID}}(\Pi^{0}_{1})\vdash\vartheta(\Omega_{k}(\ell))\in W$ for {\rm each} $\ell<\omega$.
\eprp
\bprf
Let $\beta\in M:\Leftrightarrow (K(\beta)\subset W)$.
For a formula $C(\beta)$, let $(M\to C)(\beta):\Leftrightarrow (\beta\in M\to C(\beta))$, and
${\sf J}[C](\zeta):\Leftrightarrow \forall\alpha(M\cap\alpha\subset C \to M\cap(\alpha+\Omega^{\zeta})\subset C)$.
Then we claim that
$\Pi^{0}_{k}(\Omega)\mbox{{\rm -ID}}(\Pi^{0}_{1})$ proves
$Prog(M\to C)\to Prog(M\to {\sf J}[C])$ for $\Pi^{0}_{k}(\Omega)$-formula $C$.
Argue in $\Pi^{0}_{k}(\Omega)\mbox{{\rm -ID}}(\Pi^{0}_{1})$, and assume $Prog(M\to C)$, $\forall\xi\in M\cap\zeta\,{\sf J}[C](\xi)$,
$M\cap\alpha\subset C$, $\zeta\in M$, and $\beta\in M\cap(\alpha+\Omega^{\zeta})$.
We need to show $C(\beta)$. We can assume $\beta\geq\alpha$ by $M\cap\alpha\subset C$.
When $\zeta=0$, we have $M\ni\beta=\alpha$.
Then $Prog(M\to C)$ together with $M\cap\alpha\subset C$ yields $C(\alpha)$.
Next consider the case when $\zeta$ is a limit number.
Then $\beta<\alpha+\Omega^{\xi}$ for a $\xi\in M\cap\zeta$.
$\forall\xi\in M\cap\zeta\,{\sf J}[C](\xi)$ yields $C(\beta)$.

Finally let $\zeta=\xi+1$.
Then $\xi\in M$, and we see from $\beta\in M$ that
there exists a $\gamma_{1}\in W\cap\Omega$ such that $\beta<\alpha+\Omega^{\xi}\gamma_{1}$.
We claim that $Prog(\sigma_{1})$ for
$\sigma_{1}(\gamma) : \Leftrightarrow  (M\cap(\alpha+\Omega^{\xi}\gamma)\subset C)$.
Assuming  $\forall\gamma<\gamma_{0}\, \sigma_{1}(\gamma)$, we need to show $M\cap(\alpha+\Omega^{\xi}\gamma_{0})\subset C$.
The case $\gamma_{0}=0$ follows from $M\cap\alpha\subset C$, and
the case when $\gamma_{0}$ is a limit number is readily seen.
Let $\gamma_{0}=\gamma+1$.
From ${\sf J}[C](\xi)$ and $\sigma_{1}(\gamma)$, i.e., $M\cap(\alpha+\Omega^{\xi}\gamma)\subset C$, we see that
$M\cap(\alpha+\Omega^{\xi}\gamma+\Omega^{\xi})\subset C$.
Thus $Prog(\sigma_{1})$ is shown.
Since $\sigma_{1}$ is a $\Pi^{0}_{k}(\Omega)$-formula for $k\geq 1$, we obtain
$W\subset\sigma_{1}$, i.e., $\forall\gamma\in W(M\cap(\alpha+\Omega^{\xi}\gamma)\subset C)$.
We conclude $C(\beta)$ from $\beta\in M\cap(\alpha+\Omega^{\xi}\gamma_{1})$ and $\gamma_{1}\in W$.

Thus we have shown the claim $Prog(M\to C)\to Prog(M\to{\sf J}[C])$ for $\Pi^{0}_{k}(\Omega)$-formula $C$.
Now let $C_{0}(\beta):\Leftrightarrow(\vartheta(\beta)\in W)$, and
$C_{i+1}:\equiv {\sf J}[C_{i}]$.
Then $C_{i}$ is a $\Pi^{0}_{i+1}(\Omega)$-formula.
It is clear that $Prog(M\to C_{0})$, i.e., $\forall\alpha(M\cap\alpha\subset C_{0}\to \alpha\in M\to \vartheta(\alpha)\in W)$, 
cf.\,the proof of Lemma \ref{lem:wfprf}.
By metainduction on $i\leq k$ we obtain $Prog(M\to C_{i})$.
In particular for each $\ell<\omega$,
$C_{k}(\ell)$ follows from $Prog(M\to C_{k})$.
By metainduction on $i\leq k$ we see from this that $M\cap\Omega_{i}(\ell)\subset C_{k-i}$, and hence $C_{k-i}(\Omega_{i}(\ell))$.
Therefore $C_{0}(\Omega_{k}(\ell))$, i.e., $\vartheta(\Omega_{k}(\ell))\in W$ for each $\ell<\omega$.
\eprf

The next lemma shows the easy half in Theorem \ref{th:main}.\ref{th:maina},
and the power of disjunctions of negative and positive formulas, i.e., implications of
positive formulas in the axiom (\ref{eq:Rax2}).
Note that our proof of the lemma is formalizable in $\mbox{{\rm RCA}}_{0}^{*}+(\Pi^{1}_{1}(\Pi^{0}_{3})\mbox{{\rm -CA}})^{-}$.

\blem\label{lem:Pi0WF}
$(Acc,\mbox{{\rm N}}\lor\mbox{{\rm P}})\mbox{{\rm -ID}}({\sf Acc})\vdash \vartheta(\Omega^{\ell})\in W$ for {\rm each} $\ell<\omega$.
\elem
\bprf
Argue in $(Acc,\mbox{{\rm N}}\lor\mbox{{\rm P}})\mbox{{\rm -ID}}({\sf Acc})$.
We claim that $Prog(\omega\to{\sf J}[C_{0}])$ for $C_{0}(\beta)\Leftrightarrow(\vartheta(\beta)\in W)$ and
${\sf J}[C_{0}](\ell)\Leftrightarrow \forall\alpha(M\cap\alpha\subset C_{0} \to M\cap(\alpha+\Omega^{\ell})\subset C_{0})$
with $\beta\in M\Leftrightarrow(K(\beta)\subset W)$.
${\sf J}[C_{0}](0)$ is seen from $Prog(M\subset C_{0})$.
Assuming ${\sf J}[C_{0}](\ell)$, and
$M\cap\alpha\subset C_{0}$,
we need to show $M\cap(\alpha+\Omega^{\ell+1})\subset C_{0}$.

Let $<_{lx}$ denote the lexicographic ordering on
$OT^{\prime}(\vartheta)\times OT^{\prime}(\vartheta)$, in which the first components are ordered in the ordering $<$ on $OT^{\prime}(\vartheta)$
and the second components are ordered in the $\omega$-ordering $<^{\mathbb{N}}$ on $OT^{\prime}(\vartheta)\subset\mathbb{N}$:
\[
(\xi,\gamma)<_{lx}(\zeta,\beta):\Leftrightarrow (\xi<\zeta)\lor (\xi=\zeta\land \gamma<^{\mathbb{N}}\beta)
\]
Let $W_{lx}$ denote the accessible part of $<_{lx}$, which is the least fixed point of the operator
$\forall(\xi,\gamma)<_{lx}(\zeta,\beta)\, X(\xi,\gamma)$.
Let $Prog_{lx}(X):\Leftrightarrow \forall(\zeta,\beta)[\forall(\xi,\gamma)<_{lx}(\zeta,\beta)\, X(\xi,\gamma) \to X(\zeta,\beta)]$.

\begin{equation}\label{eq:WWlx}
\zeta\in W \to \forall \beta[ (\zeta,\beta)\in W_{lx}]
\end{equation}
This follows from $Prog(D)$ for $D(\zeta)\Leftrightarrow \forall \beta[ (\zeta,\beta)\in W_{lx}]$ with the positive formula 
$D\in\mbox{Pos}\subset\mbox{{\rm N}}\lor\mbox{{\rm P}}$.
$Prog(D)$ is seen from $Acc$-induction on $\beta$.

Now let 
$\sigma_{0}(\zeta,\beta) :\Leftrightarrow  (\beta\in M \land \beta<\alpha+\Omega^{\ell}\zeta \to \vartheta(\beta)\in W)$.
We claim that $\sigma_{0}$ is progressive with respect to the lexicographic ordering $<_{lx}$, $Prog_{lx}(\sigma_{0})$.
Suppose $\forall(\xi,\gamma)<_{lx}(\zeta,\beta)\,\sigma_{0}(\xi,\gamma)$, 
$\beta\in M$ and $\beta<\alpha+\Omega^{\ell}\zeta$. 
We need to show $\vartheta(\beta)\in W$.
We can assume that $\beta=\alpha+\Omega^{\ell}\xi+\delta$ with $\xi<\zeta$ and $\delta<\Omega^{\ell}$.
We claim that $ M\cap\alpha_{0}\subset C_{0}$ for $\alpha_{0}=\alpha+\Omega^{\ell}\xi$.
Let $\gamma\in M\cap\alpha_{0}$. 
We have $(\xi,\gamma)<_{lx}(\zeta,\beta)$, $\gamma\in M$ and $\gamma<\alpha+\Omega^{\ell}\xi$.
$\sigma_{0}(\xi,\gamma)$ yields $\vartheta(\gamma)\in W$, i.e., $\gamma\in C_{0}$.
${\sf J}[C_{0}](\ell)$ yields $\beta\in M\cap(\alpha_{0}+\Omega^{\ell})\subset C_{0}$ from $M\cap\alpha_{0}\subset C_{0}$.
Thus $\vartheta(\beta)\in W$.

From $Prog_{lx}(\sigma_{0})$ we obtain $\forall(\zeta,\beta)\in W_{lx}\,\sigma_{0}(\zeta,\beta)$ for 
$\sigma_{0}\in\mbox{{\rm N}}\lor\mbox{{\rm P}}$.
By (\ref{eq:WWlx}) we conclude $\forall\zeta\in W\forall\beta\,\sigma_{0}(\zeta,\beta)$, and hence 
$M\cap(\alpha+\Omega^{\ell+1})\subset C_{0}$.

We have shown
$Prog(\omega\to{\sf J}[C_{0}])$.
By meta induction on $\ell$ we obtain ${\sf J}[C_{0}](\ell)$, and
$\vartheta(\Omega^{\ell})\in W$.
\eprf

Lemma \ref{lem:Pi0WF} shows that
\[
\vartheta(\Omega^{\omega})\leq|(Acc,\mbox{{\rm N}}\lor\mbox{{\rm P}})\mbox{{\rm -ID}}({\sf Acc})|
\leq|(Acc,\Pi^{0}_{0}(\mbox{{\rm P}}))\mbox{{\rm -ID}}({\sf Acc})|\leq|\Pi^{0}_{1}(\mbox{{\rm P}})\mbox{{\rm -ID}}({\sf Acc})|
.\]

The non-trivial halves of Theorem \ref{th:main} follow from the following theorem.
For a positive operator $\varphi(X,x)$ and a number $n$ in the least fixed point $I_{\varphi}$ of the monotonic 
operator $\omega\supset\mathcal{X}\mapsto \{n: \mathbb{N}\models\varphi[\mathcal{X},n]\}$,
$|n|_{\varphi}:=\min\{\alpha: n\in I_{\varphi}^{\alpha+1}\}$ denotes the inductive norm of $n$.
$Th(\mathbb{N})$ denotes the set of true arithmetic sentences.

\bth\label{th:main1}
\begin{enumerate}
\item\label{th:main1PN}
 For each $k\geq 0$ and positive operator $\varphi(X,x)$,
 \[
 Th(\mathbb{N})+(\Pi^{0}_{k}(\mbox{{\rm P}}),\mbox{{\rm P}}\cup\mbox{{\rm N}})\mbox{{\rm -ID}}\vdash R_{\varphi}(n)
\Rightarrow |n|_{\varphi}<\vartheta(\Omega\cdot\omega_{1+k}).
\]
\item\label{th:main1PandN}
 For each $k\geq 0$ and ${\sf Acc}$-operator $\varphi(X,x)$,
 \[
 Th(\mathbb{N})+(\Pi^{0}_{k}(\mbox{{\rm P}}),\mbox{{\rm P}}\land\mbox{{\rm N}})\mbox{{\rm -ID}}({\sf Acc})\vdash R_{\varphi}(n)
\Rightarrow |n|_{\varphi}<\vartheta(\Omega\cdot\omega_{1+k}).
\]
\item\label{th:main11}
 For each ${\sf Acc}$-operator $\varphi(X,x)$,
 \[
 Th(\mathbb{N})+\Pi^{0}_{1}(\mbox{{\rm P}})\mbox{{\rm -ID}}({\sf Acc})\vdash R_{\varphi}(n)
\Rightarrow |n|_{\varphi}<\vartheta(\Omega^{\omega}).
\]
\end{enumerate}
\end{theorem}

Our proof of Theorem \ref{th:main1} is based on an analysis through the operator controlled derivations
 due to W. Buchholz\cite{3}.
An ordinal notation system with the $\psi$-function also due to W. Buchholz\cite{2}
 (but without the exponential function below $\Omega$)
is convenient for our proof.

\bdf\label{df:Cpsi}
{\rm Let $\Omega$ be the least uncountable ordinal $\omega_{1}$, and $\varepsilon_{\Omega+1}$ the next epsilon number above $\Omega$.
Define simultaneously on ordinals $\alpha<\varepsilon_{\Omega+1}$,
operators $\mathcal{H}_{\alpha}$ on the power set of  $\varepsilon_{\Omega+1}$, and ordinals $\psi\alpha$ as follows.
Let $X\subset \varepsilon_{\Omega+1}$.
\begin{enumerate}
\item
$\{0,\Omega\}\cup X\subset \mathcal{H}_{\alpha}(X)$.
\item
If $\Omega<\beta\in\mathcal{H}_{\alpha}(X)$, then $\omega^{\beta}\in\mathcal{H}_{\alpha}(X)$.
\item
$\{\beta,\gamma\}\subset\mathcal{H}_{\alpha}(X)\Rightarrow \beta+\gamma\in\mathcal{H}_{\alpha}(X)$.
\item
$\beta\in\mathcal{H}_{\alpha}(X)\cap\alpha\Rightarrow \psi\beta\in\mathcal{H}_{\alpha}(X)$.
\end{enumerate}
Let
\[
\psi\alpha:=\min\{\beta\leq\Omega: \mathcal{H}_{\alpha}(\beta)\cap\Omega\subset\beta\}
.\]
}
\edf



It is well known that $\mathcal{H}_{\varepsilon_{\Omega+1}}(0)$ is a computable notation system, and
$\psi\alpha$ is in normal form if $G\alpha<\alpha$ for $\alpha\in\mathcal{H}_{\varepsilon_{\Omega+1}}(0)$, where
$G0=G\Omega=\emptyset$, $G(\psi\alpha)=\{\alpha\}\cup G\alpha$, $G\omega^{\alpha}=G\alpha$ and $G(\beta+\gamma)=G\beta\cup G\gamma$.
Also it is shown the following in \cite{4}.

\bprp\label{prp:Buchholz}
$\vartheta(\Omega\cdot\omega_{1+k})=\psi(\Omega^{\omega_{1+k}})$, $\vartheta(\Omega\cdot\varepsilon_{0})=\psi(\Omega^{\varepsilon_{0}})$
and $\vartheta(\Omega^{\omega})=\psi(\Omega^{\Omega^{\omega}})=\psi(\omega^{\Omega^{\omega}})$.
\eprp

Let $W$ denote the accessible part of $<$ on $\mathcal{H}_{\varepsilon_{\Omega+1}}(0)$.
The easy half in Theorem \ref{th:main}.\ref{th:mainc} follows from the following lemma.
\blem\label{lem:lowerb}
For {\rm each} $\alpha<\psi(\Omega^{\omega_{1+k}})$,
$(\Pi^{0}_{k}(\mbox{{\rm P}}),Acc)\mbox{{\rm -ID}}({\sf Acc})\vdash \alpha\in W$.
\elem
\bprf
It is clear that $Acc\mbox{{\rm -ID}}({\sf Acc})\subset(\Pi^{0}_{0}(\mbox{{\rm P}}),Acc)\mbox{{\rm -ID}}({\sf Acc})$,
and we have (a) and (b) in the proof of Lemma \ref{lem:wfprf} in hand.
The following (e) and (f) are provable in $Acc\mbox{-ID({\sf Acc})}$.
\begin{enumerate}
\item[(e)]
$G\beta<\beta\to [\forall\gamma<\beta(\mathbb{P}(\gamma)\subset W\to w(\gamma))\leftrightarrow w(\beta)]$, where
$w(\gamma):\Leftrightarrow (G\gamma<\gamma \to \psi(\gamma)\in W)$ and $\mathbb{P}(\gamma)$ denotes the set of ordinal terms $\psi\alpha$
occurring in $\gamma$.

Assume $G\beta<\beta$ and $\forall\gamma<\beta(\mathbb{P}(\gamma)\subset W\to w(\gamma))$.
By $Acc$-induction on the length of $\alpha$ we see that $\forall\alpha<\psi\beta(\alpha\in W)$.
For $\alpha=\psi\gamma$ with $G\gamma<\gamma$, $\mathbb{P}(\gamma)\subset W$ follows from the induction hypothesis and $\mathbb{P}(\gamma)<\psi\gamma$.

\item[(f)]
$Prog(E)$ for
$E(a):\Leftrightarrow (\forall\beta[\forall\gamma<\beta\, w(\gamma) \to \forall\gamma<\beta+\Omega^{a} w(\gamma)])$.

It suffices to show
$E(a+1)$ assuming $E(a)$, which follows from $Prog(D)$, and
the axiom (\ref{eq:Rax2}) for the $Acc$-formula 
$D(\zeta):\Leftrightarrow (\zeta<\Omega \to w(\beta+\Omega^{a}\zeta))$.
\end{enumerate}

From (f) we see that
$Acc\mbox{-ID({\sf Acc})}\vdash \forall\beta<\Omega^{n}\,w(\beta)$, i.e, $Acc\mbox{-ID({\sf Acc})}\vdash \forall\alpha<\psi(\Omega^{n})(\alpha\in W)$
 for each $n$.
 
In what follows argue in $(\Pi^{0}_{k}(\mbox{{\rm P}}),Acc)\mbox{-ID({\sf Acc})}$.

For a formula $A$, let
$
{\sf j}[A](\alpha) :\Leftrightarrow \forall \beta[\forall\gamma<\beta\, A(\gamma) \to \forall\gamma<\beta+\omega^{\alpha}A(\gamma)]
$.
Then $Prog(A)\to Prog({\sf j}[A])$ for $A\in\Pi^{0}_{k}(\mbox{{\rm P}})$. 

Let $E_{1}=E$ for the formula $E$ in (f), and $E_{n+1}={\sf j}[E_{n}]$.

Then $E_{n}\in\Pi^{0}_{n}(\mbox{{\rm P}})$ and $Prog(E_{k+1})$.
This yields $E_{k+1}(n)$ for each $n$.
Hence $E_{k+1-m}(\omega_{m}(n))$ for each $n$ and $m\leq k$, where
$\omega_{0}(n)=n$ and $\omega_{m+1}(n)=\omega^{\omega_{m}(n)}$, i.e., $\omega_{m}=\omega_{m}(1)$.
In particular
$E_{1}(\omega_{k}(n))$ for each $n$.
Therefore $w(\Omega^{\omega_{k}(n)})$ for each $n$.
We conclude $\forall\alpha<\psi(\Omega^{\omega_{k}(n)})(\alpha\in W)$ in $(\Pi^{0}_{k}(\mbox{{\rm P}}),Acc)\mbox{-ID({\sf Acc})}$.
\eprf

\section{Sequent calculi for weak fragments}\label{sec:finitarycal}
To establish upper bounds in Theorem \ref{th:main},
let us reformulate $Th(\mathbb{N})+(\Pi^{0}_{k}(\mbox{{\rm P}}),\mbox{{\rm P}}\cup\mbox{{\rm N}})\mbox{{\rm -ID}}$,
$Th(\mathbb{N})+(\Pi^{0}_{k}(\mbox{{\rm P}}),\mbox{{\rm P}}\land\mbox{{\rm N}})\mbox{{\rm -ID}}({\sf Acc})$ and
$Th(\mathbb{N})+\Pi^{0}_{1}(\mbox{{\rm P}})\mbox{{\rm -ID}}({\sf Acc})$ in one-sided sequent calculi.
We assume that for each predicate symbol $R$,
its complement or negation $\bar{R}$ is in the language.
For example, we have negations $\neq,\not<$ of the predicate constants $=,<$.
Logical connectives are $\lor,\land,\exists,\forall$.
Negations $\lnot A$ of formulas $A$ are defined recursively by de Morgan's law and elimination of double negations.
$A\to B$ denotes $\lnot A\lor B$ for formulas $A,B$.
$\lnot A$ is also denoted by $\bar{A}$.

The followings are \textit{initial sequents}.
\begin{enumerate}
\item
(logical initial sequent)
\[
\bar{L},L,\Gamma
\mbox{ where $L$ is a literal.}
\]

\item
(equality initial sequent)
\[
t\neq s,\bar{L}(t),L(s),\Gamma
\mbox{ for literals $L(x)$.}
\]

\item
(arithmetical initial sequent)
\[
A,\Gamma
\]
where $A$ is one of formulas $t=t$, a defining axiom for an elementary recursive function,
or a true arithmetical sentence in $\mathcal{L}$.

\end{enumerate}

Inference rules are 
$(cut)$, $(\exists),(\forall),(b\exists),(b\forall),(\lor),(\land), (R), (\bar{R})$, and $(ind)$.

\[
\infer[(cut)]{\Gam,\Del}
{
\Gam,\bar{C}
&
C,\Del
}
\]
where $C$ is the \textit{cut formula} of the $(cut)$.

\[
\infer[(\exi)]{\Gam}
{
A(t),\Gamma
}
\:
\infer[(\fal)]{\Gam}
{
A(a),\Gamma
}
\]
where $(\exists x\, A(x))\in\Gamma$ in $(\exi)$, and $a$ is an eigenvariable and $(\forall x\, A(x))\in\Gamma$.
\[
\infer[(b\exi)]{\Gam}
{
A(t),\Gamma
&
t<s,\Gamma
}
\:
\infer[(b\fal)]{\Gam}
{
a\not<s,A(a),\Gamma
}
\]
where $(\exists x<s\, A(x))\in\Gamma$ in $(b\exi)$, and $a$ is an eigenvariable and $(\forall x<s\, A(x))\in\Gamma$.
\[
\infer[(\lor)]{\Gam}
{
A_{i},\Gamma
}
\:
\infer[(\land)]{\Gam}
{
A_{0},\Gamma
&
A_{1},\Gamma
}
\]
for an $i=0,1$ with $(A_{0}\lor A_{1})\in\Gamma$ in $(\lor)$, and
$(A_{0}\land A_{1})\in\Gamma$.

 For each theory the inference rule for the predicates $R_{\varphi}$ is the following:
\[
\infer[(R)]{\Gam}
{
\varphi(R_{\varphi},t),\Gamma
}
\]
with $(R_{\varphi}(t))\in\Gamma$.

\begin{enumerate}
 \item
  For the theory $Th(\mathbb{N})+(\Pi^{0}_{k}(\mbox{{\rm P}}),\mbox{{\rm P}}\cup\mbox{{\rm N}})\mbox{{\rm -ID}}$, 
  the following $(\bar{R})$ is the inference rule for $\bar{R}_{\varphi}$:
\[
\infer[(\bar{R})]{\Gam}
{
\bar{\varphi}(\sigma,a),\sigma(a),\Gamma
&
\bar{\sigma}(t),\Gamma
}
\] 
with $(\bar{R}_{\varphi}(t))\in\Gamma$ and an eigenvariable $a$,
where 
$\varphi(X,x)$ is an $X$-positive formula, and $\sigma\in\mbox{P}\cup\mbox{N}$.

   \item
    For the theory $Th(\mathbb{N})+(\Pi^{0}_{k}(\mbox{{\rm P}}),\mbox{{\rm P}}\land\mbox{{\rm N}})\mbox{{\rm -ID}}({\sf Acc})$, 
    let $\sigma\equiv(\bar{D}\land C)$ for positive formulas $D,C$, and
  $\varphi(X,x)$ an ${\sf Acc}$-operator in (\ref{eq:Acc}).
  Then  the following $(\bar{R})$ is the inference rule for $\bar{R}_{\varphi}$:
  \[
  \infer[(\bar{R}]{\Gam}
  {
  \lnot\varphi(\bar{D},a)\lor \lnot\varphi(C,a),\sigma(a),\Gamma
  &
 \bar{\sigma}(t),\Gamma
  }
  \]
with $(\bar{R}_{\varphi}(t))\in\Gamma$ and an eigenvariable $a$.
Note that $\varphi(\bar{D},a)\land\varphi(C,a)$ is logically equivalent to $\varphi(\sigma,a)$.
  
  \item
  For the theory $Th(\mathbb{N})+\Pi^{0}_{1}(\mbox{{\rm P}})\mbox{{\rm -ID}}({\sf Acc})$, 
  let $\sigma(u)\equiv(\forall z\,\sigma_{0}(z,u))$ for $\sigma_{0}\in\Pi^{0}_{0}(\mbox{{\rm P}})$, and
  $\varphi(X,x)$ an ${\sf Acc}$-operator $\forall y\{\theta_{0}(x,y)\to t_{0}(x,y)\in X\}$ 
  with an arithmetic bounded formula $\theta_{0}(x,y)$ and a term $t_{0}(x,y)$.
  Let 
  \[
  \varphi_{\sigma}(x):\equiv[\forall w\{\theta_{0}(x, p_{0}(w))\to\sigma_{0}(p_{1}(w),t_{1}))\}]
  \]
  for $t_{1}\equiv(t_{0}(x,p_{0}(w)))$ and inverses $p_{0},p_{1}$ of a surjective pairing function.
  Note that $\varphi_{\sigma}(x)\leftrightarrow \varphi(\sigma,x)$ over {\sf EA}.
Then  the following $(\bar{R})$ is the inference rule for $\bar{R}_{\varphi}$
with $(\bar{R}_{\varphi}(t))\in\Gamma$ and an eigenvariable $a$:
\[
\infer[(\bar{R})]{\Gam}
{
\lnot\varphi_{\sigma}(a),\sigma(a),\Gamma
&
\bar{\sigma}(t),\Gamma
}
\]

  \end{enumerate}


\[
\infer[(ind)]{\Del}
{
\Delta,\theta(0)
&
\Delta,\bar{\theta}(a),\theta(a+1)
&
\bar{\theta}(t),\Delta
}
\]
where $a$ is the eigenvariable.
 \begin{enumerate}
 \item
  The \textit{induction formula} $\theta\in\Pi^{0}_{k}(\mbox{{\rm P}})$ for
  $Th(\mathbb{N})+(\Pi^{0}_{k}(\mbox{{\rm P}}),\mbox{{\rm P}}\cup\mbox{{\rm N}})\mbox{{\rm -ID}}$ and
  for $Th(\mathbb{N})+(\Pi^{0}_{k}(\mbox{{\rm P}}),\mbox{{\rm P}}\land\mbox{{\rm N}})\mbox{{\rm -ID}}({\sf Acc})$.
  \item
  $\theta\in\Pi^{0}_{1}(\mbox{{\rm P}})$
  for $Th(\mathbb{N})+\Pi^{0}_{1}(\mbox{{\rm P}})\mbox{{\rm -ID}}({\sf Acc})$.
  \end{enumerate}

Note that we can assume that when $k=0$,
$\theta\in\Pi^{0}_{0}(\mbox{{\rm P}})$ is either a formula $\exists y<t\forall z<s\bigwedge_{i}(C_{i}\to D_{i})$ 
for some positive formulas $C_{i},D_{i}$,
or its complement $\forall y<t\exists z<s\bigvee_{i}(C_{i}\land\bar{D}_{i})$.
When $k>0$, we can assume that
$\theta\in\Pi^{0}_{k}(\mbox{{\rm P}})$ is of the form $\forall x_{k} \exists x_{k-1}\cdots Q x_{1}\, \theta_{0}$,
where $Q=\forall$ if $k$ is odd, and $Q=\exists$ else, and
$\theta_{0}\in\Pi^{0}_{0}(\mbox{{\rm P}})$ is one of formulas $\exists y<t\forall z<s\bigwedge_{i}(C_{i}\to D_{i})$ and
$\forall y<t\exists z<s\bigvee_{i}(C_{i}\land\bar{D}_{i})$.

A \textit{proof} is defined from these initial sequents and inference rules.

\section{Infinitary derivations}\label{sec:infderiv}
In what follows we assume that each formula has no free variable, and 
a closed term $t$ is identified with the numeral $n$ of the value of $t$.
Furthermore assume that there occurs no bounded quantifiers in any formula.
Each bounded quantifier $\exists x<n\, B(x), \forall x<n\,B(x)$ is replaced by $\bigvee_{i<n}B(i), \bigwedge_{i<n}B(i)$, resp.
In other words, $\bigvee_{i<n}B(i), \bigwedge_{i<n}B(i)$ are formulas for formulas $\{B_{i}\}_{i<n}$.

\subsection{$\omega$-rule}
Finitary proof in the sequent calculus for $Th(\mathbb{N})+(\Pi^{0}_{k}(\mbox{{\rm P}}),\mbox{{\rm P}}\cup\mbox{{\rm N}})\mbox{{\rm -ID}}$
or for $Th(\mathbb{N})+(\Pi^{0}_{k}(\mbox{{\rm P}}),\mbox{{\rm P}}\land\mbox{{\rm N}})\mbox{{\rm -ID}}({\sf Acc})$ with $k>0$
is embedded in infinitary derivations with the $\omega$-rule:
\[
\infer{\Gam}
{
\{\Gamma,A(n):n\in\mathbb{N}\}
}
\]
with $(\forall x\, A)\in\Gamma$.

Let $I_{\varphi}^{<\Omega}:\equiv R_{\varphi}$ and $\bar{I}^{<\Omega}_{\varphi}:\equiv\bar{R}_{\varphi}$.
A formula is said to be \textit{positive} [\textit{negative}] if the predicates $\bar{I}_{\varphi}^{<\Omega}$
[the predicates $I_{\varphi}^{<\Omega}$] do not occur in it.

Definition \ref{df:SigPi} is modified as follows.
\bdf\label{df:SigPiome}
{\rm
\begin{enumerate}
\item
$\Pi^{0}_{0}(\mbox{{\rm P}})=\Sigma^{0}_{0}(\mbox{{\rm P}})$ denotes a class of formulas of the form 
$\bigvee_{i}\bigwedge_{j}(C_{ij}\to D_{ij})$ 
for some positive formulas $C_{ij},D_{ij}$,
or its complement $\bigwedge_{i}\bigvee_{j}(C_{ij}\land\bar{D}_{ij})$.

\item
If $A\in\Sigma^{0}_{k}(\mbox{{\rm P}})$, then $(\forall x\,A)\in\Pi^{0}_{k+1}(\mbox{{\rm P}})$.
If $A\in\Pi^{0}_{k}(\mbox{{\rm P}})$, then $(\exists x\,A)\in\Sigma^{0}_{k+1}(\mbox{{\rm P}})$.
\end{enumerate}
}
\edf

\bdf\label{df:dg}
{\rm
The \textit{degree} $\dg(A)<\omega$ of the formula $A\in\bigcup_{k<\omega}(\Sigma^{0}_{k}(\mbox{{\rm P}})\cup\Pi^{0}_{k}(\mbox{{\rm P}}))$ 
 is defined as follows.
\begin{enumerate}
\item
$\dg(A)=0$ if no predicate $I_{\varphi}^{<\Omega},\bar{I}_{\varphi}^{<\Omega}$ occurs in $A$.

\item
$\dg(A)=1+\min\{k: A\in\Sigma^{0}_{k}(\mbox{{\rm P}})\cup\Pi^{0}_{k}(\mbox{{\rm P}})\}$
if one of the predicates $I_{\varphi}^{<\Omega},\bar{I}_{\varphi}^{<\Omega}$ occurs in $A$.
\end{enumerate}
}
\edf

\bdf
{\rm
For finite sets $\Gamma$ of formulas, ordinals $a<\varepsilon_{0}$ and $d<\omega$,
\[
\vdash^{a}_{d}\Gamma
\]
 designates that there exists an infinitary derivation with its ordinal depth$\leq a$ and
its cut degree$<d$, where
an infinitary derivation is a well founded tree of sequents locally correct with inference rules in the
sequent calculus for $Th(\mathbb{N})+(\Pi^{0}_{k}(\mbox{{\rm P}}),\mbox{{\rm P}}\cup\mbox{{\rm N}})\mbox{{\rm -ID}}$ 
or for $Th(\mathbb{N})+(\Pi^{0}_{k}(\mbox{{\rm P}}),\mbox{{\rm P}}\land\mbox{{\rm N}})\mbox{{\rm -ID}}({\sf Acc})$ 
except
the inference rule $(\forall)$ is replaced by the $\omega$-rule.
By an infinitary derivation with cut degree$<d$, we mean a derivation in which $\dg(A)<d$ for every cut formula $A$.
}
\edf

Let $\Gamma[\vec{a}]$ be a sequent in the language of $\mathcal{L}(\mbox{{\rm ID}})$, where
$\vec{a}=(a_{1},\ldots,a_{p})$ is a list of free variables occurring in the sequent $\Gamma[\vec{a}]$.
For lists $\vec{n}=(n_{1},\ldots,n_{p})\subset \mathbb{N}$ of natural numbers, 
$\Gamma^{*}[\vec{n}]=\{A^{*}[\vec{n}]:A\in\Gamma\}$ and 
$A^{*}[\vec{n}]$ denotes the result of replacing
every occurrence of the variable $a_{i}$ in the list $\vec{a}$ by the natural number $n_{i}$, 
and
every occurrence of bounded quantifies $\exists x<n\, B(x), \forall x<n\,B(x)$
by $\bigvee_{i<n}B^{*}(i), \bigwedge_{i<n}B^{*}(i)$, resp.

\blem\label{lem:preembed}{\rm (Pre-embedding)}
\begin{enumerate}
\item
If $Th(\mathbb{N})+(\Pi^{0}_{k}(\mbox{{\rm P}}),\mbox{{\rm P}}\cup\mbox{{\rm N}})\mbox{{\rm -ID}}\vdash \Gamma[\vec{a}]$ for $k>0$, 
then there exists $a<\omega_{1+k}$ such that $\vdash_{2}^{a}\Gamma^{*}[\vec{n}]$ for any $\vec{n}$.
\item
If $Th(\mathbb{N})+(\Pi^{0}_{k}(\mbox{{\rm P}}),\mbox{{\rm P}}\land\mbox{{\rm N}})\mbox{{\rm -ID}}({\sf Acc})\vdash \Gamma[\vec{a}]$ for $k>0$, 
then there exists $a<\omega_{1+k}$ such that $\vdash_{2}^{a}\Gamma^{*}[\vec{n}]$ for any $\vec{n}$.
\end{enumerate}
\elem
\bprf
Consider $Th(\mathbb{N})+(\Pi^{0}_{k}(\mbox{{\rm P}}),\mbox{{\rm P}}\cup\mbox{{\rm N}})\mbox{{\rm -ID}}$.
Let $P$ be a proof of the sequent $\Gamma[\vec{a}]$.
By eliminating $(cut)$'s partially we may assume that any cut formula in $P$ is either an arithmetical formula or 
an atomic formulas $R_{\varphi}(t)$.
We see easily that there exists $c<\omega^{2}$ such that $\vdash_{2+k}^{c}\Gamma^{*}[\vec{n}]$ for any $\vec{n}$
since $\dg(\theta)\leq 1+k$ for the induction formula $\theta\in\Pi^{0}_{k}(\mbox{{\rm P}})$.
By cut-elimination we obtain 
$\vdash_{2}^{a}\Gamma^{*}[\vec{n}]$ for $a=2_{k}(c)<\omega_{1+k}$ with $2_{0}(c)=c$ and $2_{m+1}(c)=2^{2_{m}(c)}$.

The lemma for $Th(\mathbb{N})+(\Pi^{0}_{k}(\mbox{{\rm P}}),\mbox{{\rm P}}\land\mbox{{\rm N}})\mbox{{\rm -ID}}({\sf Acc})$ is similarly seen.
\eprf

\subsection{Operator controlled derivations}
In this subsection let us introduce operator controlled derivations, and prove the remaining halves in Theorem \ref{th:main}.

The language $\mathcal{L}^{\infty}(\mbox{ID})$ for the next infinitary calculus
 is obtained from the language $\mathcal{L}(\mbox{ID})$ 
by deleting free variables,
and adding unary predicate symbols $I_{\varphi}^{<\alpha}, \bar{I}_{\varphi}^{<\alpha}$ for each positive operator $\varphi$
and each $\alpha<\vartheta(\Omega^{\omega})=\psi(\omega^{\Omega^{\omega}})$.
Recall that $I_{\varphi}^{<\Omega}:\equiv R_{\varphi}$.

A formula in the language $\mathcal{L}^{\infty}(\mbox{ID})$
 is said to be \textit{positive} [\textit{negative}] if the predicates $\bar{I}_{\varphi}^{<\Omega}$
[the predicates $I_{\varphi}^{<\Omega}$] do not occur in it.
In these formulas predicates $I_{\varphi}^{<\alpha}, \bar{I}_{\varphi}^{<\alpha}$ may occur.
Definition \ref{df:SigPiome} of classes $\Pi^{0}_{k}(\mbox{{\rm P}}),\Sigma^{0}_{k}(\mbox{{\rm P}})$ of formulas and
Definition \ref{df:dg} of the degree $\dg(A)$ 
of formulas $A$
are modified according to this enlargement of positive/negative formulas.
Specifically predicates $I_{\varphi}^{<\alpha}, \bar{I}_{\varphi}^{<\alpha}$ may occur in formulas $A$ with $\dg(A)=0$.

A closed term $t$ is identified with the numeral $n$ of the value of $t$.
$\Gamma,\Delta,\ldots$ denote finite sets of formulas, \textit{sequents}.

$I_{\varphi}^{<\alpha}$ is intended to denote the union $\bigcup_{\beta<\alpha}I_{\varphi}^{\beta}$ 
of the $\beta$-th stage $I_{\varphi}^{\beta}=\{n\in\mathbb{N}: \varphi(I_{\varphi}^{<\beta},n)\}$ of the least fixed point $I^{<\Omega}_{\varphi}$,
and $ \bar{I}_{\varphi}^{<\alpha}$ its complement.
Thus for any ordinal $\alpha$ and any natural number $n$
\[
I_{\varphi}^{<\alpha}(n)\leftrightarrow \exists \beta<\alpha\, \varphi(I_{\varphi}^{<\beta},n)
.\]


For a sequent $\Gamma$, 
${\sf k}(\Gamma)$ denotes the set of ordinals $\alpha<\Omega$ such that one of predicates $I_{\varphi}^{<\alpha}, \bar{I}_{\varphi}^{<\alpha}$
occurs in a formula in the set $\Gamma$.
$\mathcal{H}[\Theta](X):=\mathcal{H}(\Theta\cup X)$ for sets $\Theta, X$ of ordinals and operators $\mathcal{H}:X\mapsto\mathcal{H}(X)$
on the sets $X$ of ordinals.

\bdf\label{df:controlder}
{\rm
\textit{Inductive definition} of $\mathcal{H}
\vdash^{a}_{d}\Gamma$.

Let $\Gamma$ be a sequent, $a<\Omega\cdot\varepsilon_{0}$ and $d\leq 3$.
$\mathcal{H}\vdash^{a}_{d}\Gamma$ holds if
\begin{equation}\label{eq:controlder}
\{a\}\cup{\sf k}(\Gamma)\subset\mathcal{H}
\end{equation}
and one of the followings holds:
\begin{description}
 \item[(initial)]
 
There exists a true arithmetic formula $A\in\mathcal{L}$ in $\Gamma$.

\item[$(\bigvee)$]
There exist a formula $(\bigvee_{i<n}A_{i})\in\Gamma$ with $n>0$, $a_{0}<a$ and $i<n$ such that
$\mathcal{H}\vdash^{a_{0}}_{d}\Gamma,A_{i}$.

\item[$(\bigwedge)$]
There exist a formula $(\bigwedge_{i<n}A_{i})\in\Gamma$ and an $a_{0}<a$ such that
$\forall i<n\{\mathcal{H}\vdash^{a_{0}}_{d}\Gamma,A_{i}\}$.

\item[$(\exists)$]
There exist a formula $(\exists x A(x))\in\Gamma$, $n\in\omega$ and $a_{0}<a$ such that
$\mathcal{H}\vdash^{a_{0}}_{d}\Gamma,A(n)$.

\item[$(\forall^{\omega})$]
There exist a formula $(\forall x A(x))\in\Gamma$ and a sequence of ordinals $\{a_{n}\}_{n\in\mathbb{N}}$ such that $\forall n(a_{n}<a)$ and
$\forall n\{\mathcal{H}\vdash^{a_{n}}_{d}\Gamma,A(n)\}$.

\item[$(I^{<})$]
There exist $\alpha\leq\Omega$, $(I_{\varphi}^{<\alpha}(n))\in\Gamma$, $\beta<\alpha$ and
 $a_{0}<a$,
\begin{equation}\label{eq:controlder1}
\mbox{if } X \mbox{ occurs in } \varphi(X,n), \mbox{ then } \beta<a
\end{equation}
and
$\mathcal{H}^{\prime}\vdash^{a_{0}}_{d}\Gamma, \varphi(I_{\varphi}^{<\beta},n)$,
where $\mathcal{H}^{\prime}=\mathcal{H}[\{\beta\}]$ if $X$ occurs in $\varphi(X,n)$, and
$\mathcal{H}^{\prime}=\mathcal{H}$ else.

\item[$(\bar{I}^{<})$]
There exist $\alpha\leq\Omega$, $(\bar{I}_{\varphi}^{<\alpha}(n))\in\Gamma$ and 
a sequence of ordinals $\{a_{\beta}\}_{\beta<\alpha}$ such that
$\forall\beta<\alpha(a_{\beta}<a)$ and
$\forall\beta<\alpha\{\mathcal{H}[\{\beta\}]\vdash^{a_{\beta}}_{d}\Gamma,\lnot \varphi(I_{\varphi}^{<\beta},n)\}$.

\item[$(Cl)$]
There exist a formula $(n\in I_{\varphi}^{<\Omega})\in\Gamma$ and $a_{0}<a$ such that
$\mathcal{H}\vdash^{a_{0}}_{d}\Gamma,\varphi(I_{\varphi}^{<\Omega},n)$.

\item[$(cut)$]
There exist a formula $C$ and $a_{0}<a$ such that $\dg(C)<d$,
$\mathcal{H}\vdash^{a_{0}}_{d}\Gamma,\lnot C$ and
$\mathcal{H}\vdash^{a_{0}}_{d}C,\Gamma$.
\end{description}
}
\edf

\blem\label{th:embedregthm1}{\rm (Embedding 1)}\\
If $Th(\mathbb{N})+(\Pi^{0}_{k}(\mbox{{\rm P}}),\mbox{{\rm P}}\cup\mbox{{\rm N}})\mbox{{\rm -ID}}\vdash \Gamma[\vec{a}]$, 
there exist an $a<\Omega\cdot\omega_{1+k}$ such that 
for any $\vec{n}\subset \mathbb{N}$ and any operator $\mathcal{H}=\mathcal{H}_{\gamma}$ defined in Definition \ref{df:Cpsi} 
with $\gamma\geq 2$,
 $\mathcal{H}\vdash_{2}^{a}\Gamma^{*}[\vec{n}]$.
\elem
\bprf
Note that $1=\psi 0,\omega=\psi 1\in\mathcal{H}_{2}$.

First consider the case $k=0$.
Pick a finitary proof of the sequent $\Gamma[\vec{a}]$ in 
$Th(\mathbb{N})+(\Pi^{0}_{0}(\mbox{{\rm P}}),\mbox{{\rm P}}\cup\mbox{{\rm N}})\mbox{{\rm -ID}}$.
We show the lemma by induction on the depth of the finitary proof.

Logical initial sequents $\Gamma,\bar{R}_{\varphi}(t),R_{\varphi}(t)$ turns to
$\mathcal{H}\vdash^{\Omega}_{0}\Gamma^{*},\bar{I}_{\varphi}^{<\Omega}(n),I_{\varphi}^{<\Omega}(n)$, which in turn follows from
$\mathcal{H}[\{\alpha\}]\vdash^{f(\alpha)}_{0}\bar{I}_{\varphi}^{<\alpha}(n),I_{\varphi}^{<\alpha}(n)$ for any $\alpha<\Omega$ and any $n\in\mathbb{N}$
with $f(\alpha)=k\alpha$ for a $k<\omega$.
\[
\infer[(ind)]{\Del}
{
\Delta,\theta(0)
&
\Delta,\bar{\theta}(a),\theta(a+1)
&
\bar{\theta}(t),\Delta
}
\]

We can assume that the bounded formula $\theta(a)$ is of the form $\exists x<t\forall y<s\bigwedge_{k<m}(\bar{C}_{k}\land D_{k})$
for positive formulas $C_{k},D_{k}$.
Then $\theta^{*}(i)\equiv\bigvee _{j<p}\bigwedge_{k<q}(\bar{C}_{ijk}\land D_{ijk})$ with $\dg(\theta^{*}(i))=1$.
The inference $(ind)$
turns to a series of $(cut)$'s of cut formulas $\theta^{*}(i)$ for an $a<\Omega\cdot\omega$. 
From $\mathcal{H}\vdash_{2}^{a}\Delta^{*},\theta^{*}(0)$ and
$\mathcal{H}\vdash_{2}^{a}\Delta^{*},\bar{\theta}^{*}(0),\theta^{*}(1)$, infer
$\mathcal{H}\vdash_{2}^{a+1}\Delta^{*},\theta^{*}(1)$, and so on.
From $\mathcal{H}\vdash_{2}^{a+n}\Delta^{*},\theta^{*}(n)$ and
$\mathcal{H}\vdash_{2}^{a}\bar{\theta}^{*}(n),\Delta^{*}$, infer
$\mathcal{H}\vdash_{2}^{a+\omega}\Delta^{*}$.

Next consider
\[
\infer[(R)]{\Gam}
{
\varphi(R_{\varphi},t),\Gamma
}
\]

From $\mathcal{H}\vdash_{2}^{a}\varphi(I^{<\Omega}_{\varphi},n),\Gamma$, infer
$\mathcal{H}\vdash_{2}^{a+1}\Gamma$ by $(Cl)$.

Finally consider
\[
\infer[(\bar{R})]{\Gam}
{
\bar{\varphi}(\sigma,a),\sigma(a),\Gamma
&
\bar{\sigma}(t),\Gamma
}
\] 
where $\bar{R}_{\varphi}(t)\in\Gamma$ and $\sigma\in\mbox{P}\cup\mbox{N}$.

By the induction hypothesis we have a $c<\Omega\cdot\omega$ such that
$\mathcal{H}\vdash_{2}^{c} \bar{\varphi}(\sigma^{*},n), \sigma^{*}(n),\Gamma^{*}$
and $\mathcal{H}\vdash_{2}^{c} \bar{\sigma}^{*}(n),\Gamma^{*}$.
We show by induction on $\alpha<\Omega$ that for $f(\alpha)=c+\omega\alpha+1$ and any $n\in\mathbb{N}$
\begin{equation}\label{eq:trind}
\mathcal{H}[\{\alpha\}]\vdash^{f(\alpha)}_{2}\Gamma^{*},\bar{I}_{\varphi}^{<\alpha}(n),\sigma^{*}(n)
\end{equation}
By the induction hypothesis we have
$\mathcal{H}[\{\beta\}]\vdash^{f(\beta)}_{2}\Gamma^{*},\bar{I}_{\varphi}^{<\beta}(n),\sigma^{*}(n)$
for any $\beta<\alpha$ and $n\in\mathbb{N}$.
From this we see that $\mathcal{H}[\{\beta\}]\vdash^{f(\beta)+m}_{2}\Gamma^{*},\bar{\varphi}(I^{<\beta},n),\varphi(\sigma^{*},n)$
for some $m<\omega$.
$(\bar{I}^{<})$ yields $\mathcal{H}[\{\alpha\}]\vdash^{c+\omega\alpha}_{2}\Gamma^{*},\bar{I}_{\varphi}^{<\alpha}(n),\varphi(\sigma^{*},n)$.
A $(cut)$ with $\mathcal{H}\vdash_{2}^{c} \bar{\varphi}(\sigma^{*},n), \sigma^{*}(n),\Gamma^{*}$
yields
$\mathcal{H}[\{\alpha\}]\vdash^{f(\alpha)}_{2}\Gamma^{*}, \bar{I}_{\varphi}^{<\alpha}(n),\sigma^{*}(n)$.
Here note that $\varphi(\sigma^{*},n)\in\mbox{P}\cup\mbox{N}$ and
\begin{equation}\label{eq:sigPN}
\varphi(\sigma^{*},n)\in\Pi^{0}_{0}(\mbox{{\rm P}})
\end{equation}
with $\dg(\varphi(\sigma^{*},n))\leq 1$.

From (\ref{eq:trind}) and $(\bar{I}^{<})$ we have
$\mathcal{H}\vdash^{c+\Omega}_{2}\Gamma^{*},\bar{I}_{\varphi}^{<\Omega}(n),\sigma^{*}(n)$.
Finally a $(cut)$ with $\dg(\sigma^{*}(n))\leq 1$ yields 
$\mathcal{H}\vdash^{c+\Omega+1}_{2}\Gamma^{*}$ for $(\bar{I}_{\varphi}^{<\Omega}(n))\in\Gamma^{*}$
and $c+\Omega+1<\Omega\cdot\omega$.

Next consider the case $k>0$.
By Lemma \ref{lem:preembed}
 there exists $a<\omega_{1+k}$ such that $\vdash_{2}^{a}\Gamma^{*}[\vec{n}]$ for any $\vec{n}$.
We show by induction on $a$ that
\[
\vdash_{2}^{a}\Gamma \Rightarrow \exists c\leq\Omega\cdot (1+a)(\mathcal{H}[{\sf k}(\Gamma)]\vdash^{c}_{2}\Gamma)
\]
Consider an $\omega$-rule.
Let $(\forall x A(x))\in\Gamma$ and $\vdash^{a_{n}}_{2}\Gamma,A(n)$ with $a_{n}<a$ for any $n$.
By the induction hypothesis we have $\mathcal{H}[{\sf k}(\Gamma)]\vdash^{c_{n}}_{2}\Gamma,A(n)$ for $c_{n}\leq\Omega\cdot(1+ a_{n})$.
Then by $(\forall^{\omega})$ we obtain
$\mathcal{H}[{\sf k}(\Gamma)]\vdash^{c}_{2}\Gamma$ for $c=\sup_{n}\{\Omega\cdot(1+a_{n})+1\}\leq\Omega\cdot(1+a)$.
\eprf

\blem\label{th:embedregthm2}{\rm (Embedding 2)}\\
If $Th(\mathbb{N})+(\Pi^{0}_{k}(\mbox{{\rm P}}),\mbox{{\rm P}}\land\mbox{{\rm N}})\mbox{{\rm -ID}}({\sf Acc})\vdash \Gamma[\vec{a}]$, 
there exist an $a<\Omega\cdot\omega_{1+k}$ such that 
for any $\vec{n}\subset \mathbb{N}$ and any operator $\mathcal{H}=\mathcal{H}_{\gamma}$ defined in Definition \ref{df:Cpsi} 
with $\gamma\geq 2$,
 $\mathcal{H}\vdash_{2}^{a}\Gamma^{*}[\vec{n}]$.
\elem
\bprf
This is seen as in Lemma \ref{th:embedregthm1}.
Consider
  \[
  \infer[(\bar{R}]{\Gam}
  {
  \lnot\varphi(\bar{D},a)\lor \lnot\varphi(C,a),\sigma(a),\Gamma
  &
 \bar{\sigma}(t),\Gamma
  }
  \]
where $\bar{R}_{\varphi}(t)\in\Gamma$ and $\sigma\equiv(\bar{D}\land C)$ with positive formulas $D,C$.
As in (\ref{eq:trind}) we see for a $c<\Omega\cdot\omega_{1+k}$ and $f(\alpha)=c+\omega\alpha+1$ 
that $\mathcal{H}[\{\alpha\}]\vdash^{f(\alpha)}_{2}\Gamma^{*},\bar{I}_{\varphi}^{<\alpha}(n),\sigma^{*}(n)$ for any $n\in\mathbb{N}$.
Note that the cut formulas $\varphi(\bar{D}^{*},n)\land\varphi(C^{*},n)$ and $\sigma^{*}(n)$ arise, cf.\,(\ref{eq:sigPN}).
We have
$\varphi(\bar{D}^{*},n)\land\varphi(C^{*},n), \sigma^{*}(n)\in \Pi^{0}_{0}(\mbox{{\rm P}})$ and
$\dg(\varphi(\bar{D}^{*},n)\land\varphi(C^{*},n)),\dg(\sigma^{*}(n))\leq 1$.
\eprf

\blem\label{th:embedregthm3}{\rm (Embedding 3)}\\
If $Th(\mathbb{N})+\Pi^{0}_{1}(\mbox{{\rm P}})\mbox{{\rm -ID}}({\sf Acc})\vdash \Gamma[\vec{a}]$, 
there exist an $a<\Omega\cdot\omega$ such that 
for any $\vec{n}\subset \mathbb{N}$ and any operator $\mathcal{H}=\mathcal{H}_{\gamma}$ defined in Definition \ref{df:Cpsi}  
with $\gamma\geq 2$,
 $\mathcal{H}\vdash_{3}^{a}\Gamma^{*}[\vec{n}]$.
\elem
\bprf
This is seen as in Lemma \ref{th:embedregthm1} for $k=0$.
Note that the cut formula $\theta^{*}(i)\in\Pi^{0}_{1}(\mbox{{\rm P}})$ arises from $(ind)$ with $\dg(\theta^{*}(i))\leq 2$.
Consider
\[
\infer[(\bar{R})]{\Gam}
{
\bar{\varphi}(\sigma,a),\sigma(a),\Gamma
&
\bar{\sigma}(t),\Gamma
}
\] 
where $\bar{R}_{\varphi}(t)\in\Gamma$ and $\sigma,\varphi_{\sigma}(a)\in\Pi^{0}_{1}(\mbox{{\rm P}})$.
As in (\ref{eq:trind}) we see for a $c<\Omega\cdot\omega$ and $f(\alpha)=c+\omega\alpha+1$ 
that $\mathcal{H}[\{\alpha\}]\vdash^{f(\alpha)}_{3}\Gamma^{*},\bar{I}_{\varphi}^{<\alpha}(n),\sigma^{*}(n)$ for any $n\in\mathbb{N}$.
For the cut formulas $\varphi_{\sigma}^{*}(n)$ and $\sigma^{*}(n)$, we have
$\varphi_{\sigma}^{*}(n), \sigma^{*}(n)\in \Pi^{0}_{1}(\mbox{{\rm P}})$ and
$\dg(\varphi_{\sigma}^{*}(n)), \dg(\sigma^{*}(n))\leq 2$.
\eprf

\blem\label{lem:CE}
For any operator $\mathcal{H}=\mathcal{H}_{\gamma}$ defined in Definition \ref{df:Cpsi}  
with $\gamma\geq 2$,
if $\mathcal{H}\vdash^{a}_{3}\Gamma$, then $\mathcal{H}\vdash^{\omega^{a}}_{2}\Gamma$.
\elem

In the following lemmas $\Gamma^{(b)}=\{A^{(b)}:A\in\Gamma\}$, while $A^{(b)}$ is obtained from $A$ by replacing some
\textit{positive} occurrences of $I_{\varphi}^{<\Omega}$ by $I_{\varphi}^{<b}$.

\blem\label{lem:bounding}{\rm (Bounding)}
\\
Let $\mathcal{H}\vdash^{a}_{1}\Gamma$ for $a<\Omega$ and $\Gamma\subset\mbox{{\rm Pos}}$.
Then $\mathcal{H}\vdash^{a}_{1}\Gamma^{(b)}$ for $a\leq b\in\mathcal{H}\cap\Omega$.
\elem
\bprf
This is seen by induction on $a<\Omega$.

Suppose $\mathcal{H}\vdash^{a}_{1}\Gamma$ follows from $(I^{<})$ 
so that $(I_{\varphi}^{<\Omega}(n))\in\Gamma$ and $\mathcal{H}[\{\gamma\}]\vdash^{a_{\gamma}}_{1}\Gamma,\varphi(I_{\varphi}^{<\gamma},n)$
with $\gamma<\Omega$, $a_{\gamma}<a$ and $\gamma<a$ if $X$ occurs in $\varphi(X,n)$, (\ref{eq:controlder1}).
Then by the induction hypothesis we have $\mathcal{H}[\{\gamma\}]\vdash^{a_{\gamma}}_{1}\Gamma^{(b)},\varphi(I_{\varphi}^{<\gamma},n)$.
By $(I^{<})$ we obtain
$\mathcal{H}\vdash^{a}_{1}\Gamma^{(b)}$ for $\gamma<a\leq b$.

Suppose $\mathcal{H}\vdash^{a}_{1}\Gamma$ follows from $(Cl)$ with $(I_{\varphi}^{<\Omega}(n))\in\Gamma$ and
$\mathcal{H}\vdash^{a_{0}}_{1}\Gamma,\varphi(I_{\varphi}^{<\Omega},n)$ for an $a_{0}<a$.
By the induction hypothesis we have $\mathcal{H}\vdash^{a_{0}}_{1}\Gamma^{(b)},\varphi(I_{\varphi}^{<a_{0}},n)$ for $a_{0}<a\leq b$ and $a_{0}\in\mathcal{H}$.
An $(I^{<})$ yields $\vdash^{a}_{1}\Gamma^{(b)}$.
\eprf

\blem\label{lem:collaps}{\rm (Collapsing)}\\
Let $\gamma\in\mathcal{H}_{\gamma}$ and $\Gamma\subset\mbox{{\rm Pos}}$.
Assume $\mathcal{H}_{\gamma}\vdash_{2}^{a}\Gamma$.
Then 
$\mathcal{H}_{\hat{a}+1}\vdash^{\psi\hat{a}}_{1}\Gamma$ for $\hat{a}=\gamma+\omega^{\Omega+a}$.
\elem
\bprf
We show the lemma by induction on $a$.

First let us verify the condition (\ref{eq:controlder}) in $\mathcal{H}_{\hat{a}+1}\vdash^{\psi\hat{a}}_{1}\Gamma$.
From $\gamma<\hat{a}+1$ we see ${\sf k}(\Gamma)\subset\mathcal{H}_{\gamma}\subset\mathcal{H}_{\hat{a}+1}$.
Also by $\{\gamma,a\}\subset\mathcal{H}_{\gamma}$ we have
$\hat{a}=\gamma+\omega^{\Omega+a}\in \mathcal{H}_{\gamma}\subset\mathcal{H}_{\hat{a}}$ and
$\psi\hat{a}\in\mathcal{H}_{\hat{a}+1}$. 
From $\hat{a}\in\mathcal{H}_{\hat{a}}$ we see that if $a_{0}<a$ and
$\mathcal{H}_{\gamma}\vdash^{a_{0}}_{2}\Gamma_{0}$, then
$\psi\widehat{a_{0}}<\psi\hat{a}$.
\\

\noindent
{\bf Case 1} $(\bar{I}^{<})$:
For $a_{\beta}<a$, $\mathcal{H}_{\gamma}\vdash^{a}_{2}\Gamma, \bar{I}_{\varphi}^{<\alpha}(n)$ follows from
\\
$\{\mathcal{H}_{\gamma}[\{\beta\}]\vdash^{a_{\beta}}_{2}\Gamma, \bar{I}_{\varphi}^{<\alpha}(n),\lnot\varphi(I_{\varphi}^{<\beta},n) :\beta<\alpha\}$.
From $(\bar{I}_{\varphi}^{<\alpha}(n))\in\mbox{Pos}$ we see $\alpha<\Omega$.
We claim
\begin{equation}\label{eq:Collapsingthm}
\forall\beta<\alpha(\beta\in\mathcal{H}_{\gamma})
\end{equation} 
Let  $\beta<\alpha$.
We have $\Omega>\alpha\in{\sf k}(\bar{I}_{\varphi}^{<\alpha}(n))\subset\mathcal{H}_{\gamma}$, which yields
$\beta<\alpha\in \mathcal{H}_{\gamma}(0)\cap\Omega=\psi\gamma$, and $\beta\in\mathcal{H}_{\gamma}$.
 (\ref{eq:Collapsingthm})  yields $\mathcal{H}_{\gamma}[\{\beta\}]=\mathcal{H}_{\gamma}$.
By the induction hypothesis we obtain for any $\beta<\alpha$, 
$\mathcal{H}_{\widehat{a_{\beta}}+1}\vdash^{\psi\widehat{a_{\beta}}}_{1}\Gamma,\bar{I}_{\varphi}^{<\alpha}(n),\lnot\varphi(I_{\varphi}^{<\beta},n)$, where
$\widehat{a_{\beta}}=\gamma+\omega^{\Omega+a_{\beta}}$ and $\psi\widehat{a_{\beta}}<\psi\hat{a}$.
We conclude $\mathcal{H}_{\hat{a}+1}\vdash^{\psi\hat{a}}_{1}\Gamma,\bar{I}_{\varphi}^{<\alpha}(n)$ by $(\bar{I}^{<})$.
\\
{\bf Case 2} $(I^{<})$:
For $\beta<\min\{\alpha,a\}$ and $\alpha\leq\Omega$,
$\mathcal{H}_{\gamma}\vdash^{a}_{2}\Gamma,I_{\varphi}^{<\alpha}(n)$ follows from
$\mathcal{H}_{\gamma}\vdash^{a_{0}}_{2}\Gamma, I_{\varphi}^{<\alpha}(n), \varphi(I_{\varphi}^{<\beta},n)$.
If $X$ occurs in $\varphi(X,n)$, then by (\ref{eq:controlder}) we have
$\Omega>\beta\in{\sf k}(\varphi(I_{\varphi}^{<\beta},n))\subset\mathcal{H}_{\gamma}$, and $\beta<\psi\gamma\leq\psi\hat{a}$.
For $\widehat{a_{0}}=\gamma+\omega^{\Omega+a_{0}}$ we obtain 
$\mathcal{H}_{\widehat{a_{0}}+1}\vdash^{\psi\widehat{a_{0}}}_{1}\Gamma, I_{\varphi}^{<\alpha}(n), \varphi(I_{\varphi}^{<\beta},n)$ by the induction hypothesis.
$(I^{<})$ with $\beta<\psi\hat{a}$ for  (\ref{eq:controlder1}) yields
$\mathcal{H}_{\hat{a}+1}\vdash^{\psi\hat{a}}_{1}\Gamma,I_{\varphi}^{<\alpha}(n)$.
\\
{\bf Case 3}.
$\mathcal{H}_{\gamma}\vdash^{a}_{2}\Gamma$ follows by a $(cut)$ from
\begin{equation}\label{eq:case3.0}
\mathcal{H}_{\gamma}\vdash^{a_{0}}_{2}\Gamma, \bigwedge_{i}(C_{i}\lor \bar{D}_{i})
\end{equation}
and
\begin{equation}\label{eq:case3.1}
\mathcal{H}_{\gamma}\vdash^{a_{0}}_{2} \bigvee_{i}(\bar{C}_{i}\land D_{i}),\Gamma
\end{equation}
for $a_{0}<a$ and positive formulas $C_{i},D_{i}$.
For the sake of simplicity let us assume $i=0,1$.
In the Appendix \ref{sec:appendix} the general case is treated.

Let $b_{m}=\gamma+\omega^{\Omega+a_{0}}\cdot m$ and 
$\beta_{m}=\psi(b_{m})$ for $m=1,2,\ldots,5$.
We have $\beta_{m}<\beta_{m+1}$.
By inversion on (\ref{eq:case3.1}) we have
\[
\mathcal{H}_{\gamma}\vdash^{a_{0}}_{2} D_{0},D_{1},\Gamma
\]
By the induction hypothesis we obtain
\[
\mathcal{H}_{b_{1}+1}\vdash^{\beta_{1}}_{1} D_{0},D_{1},\Gamma
\]
From $\beta_{1}\in\mathcal{H}_{b_{1}+1}$ and Bounding lemma \ref{lem:bounding} we obtain
\begin{equation}\label{eq:case3.2}
\mathcal{H}_{b_{1}+1}\vdash^{\beta_{1}}_{1}  D_{0}^{(\beta_{1})},D_{1}^{(\beta_{1})},\Gamma
\end{equation}
For each $i=0,1$ we have by inversion on (\ref{eq:case3.0}) 
\[
\mathcal{H}_{b_{1}+1}\vdash^{a_{0}}_{2}\Gamma, C_{i}, \bar{D}_{i}^{(\beta_{1})}
\]
From $b_{1}+1\in\mathcal{H}_{b_{1}+1}$ and the induction hypothesis we obtain 
\[
\mathcal{H}_{b_{2}+1}\vdash^{\beta_{2}}_{1}\Gamma,
 C_{i}, \bar{D}_{i}^{(\beta_{1})}
\]
Once again by Bounding lemma \ref{lem:bounding} we obtain
\begin{equation}\label{eq:case3.3}
\mathcal{H}_{b_{2}+1}\vdash^{\beta_{2}}_{1}\Gamma, 
C_{i}^{(\beta_{2})}, \bar{D}_{i}^{(\beta_{1})}
\end{equation}
Again by inversion on (\ref{eq:case3.1}) we obtain
\[
\mathcal{H}_{b_{2}+1}\vdash^{a_{0}}_{2} \bar{C}^{(\beta_{2})}_{i}, D_{1-i},\Gamma
\]
and the induction hypothesis yields 
\begin{equation}\label{eq:case3.4}
\mathcal{H}_{b_{3}+1}\vdash^{\beta_{3}}_{1} \bar{C}^{(\beta_{2})}_{i},D_{1-i}^{(\beta_{3})},
\Gamma
\end{equation}

From (\ref{eq:case3.0}) we see that
\begin{equation}\label{eq:case3.5}
\mathcal{H}_{b_{4}+1}\vdash^{\beta_{4}}_{1}\Gamma, 
C_{1-i}^{(\beta_{4})},
\bar{D}_{1-i}^{(\beta_{3})}
\end{equation}
and from (\ref{eq:case3.1})  we see that
\begin{equation}\label{eq:case3.6}
\mathcal{H}_{b_{5}+1}\vdash^{\beta_{5}}_{1}\Gamma,\bar{C}_{i}^{(\beta_{2})},\bar{C}_{1-i}^{(\beta_{4})}
\end{equation}
A $(cut)$ with (\ref{eq:case3.5}) and (\ref{eq:case3.6}) yields
\begin{equation}\label{eq:case3.7}
\mathcal{H}_{b_{5}+1}\vdash^{\beta_{5}+1}_{1}\Gamma,
\bar{C}_{i}^{(\beta_{2})}, \bar{D}_{1-i}^{(\beta_{3})}
\end{equation}
Another $(cut)$ with (\ref{eq:case3.7}) and (\ref{eq:case3.4}) yields
\begin{equation}\label{eq:case3.8}
\mathcal{H}_{b_{5}+1}\vdash^{\beta_{5}+2}_{1}\Gamma,\bar{C}_{i}^{(\beta_{2})}
\end{equation}
One more $(cut)$ with (\ref{eq:case3.8}) and (\ref{eq:case3.3}) yields for each $i=0,1$
\begin{equation}\label{eq:case3.9}
\mathcal{H}_{b_{5}+1}\vdash^{\beta_{5}+3}_{1}\Gamma,\bar{D}_{i}^{(\beta_{1})}
\end{equation}
Finally several $(cut)$'s with (\ref{eq:case3.9}) and (\ref{eq:case3.2}) yields
\[
\mathcal{H}_{b_{5}+1}\vdash^{\beta_{5}+5}_{1}\Gamma
\]
Here we have
$b_{5}=\gamma+\omega^{\Omega+a_{0}}\cdot 5<\gamma+\omega^{\Omega+a}=\hat{a}$ and 
$\beta_{5}=\psi(b_{5})<\psi(\gamma+\omega^{\Omega+a})=\psi\hat{a}$.

All other cases are seen easily from the induction hypothesis.
\eprf

\noindent
({\bf Proof} of Theorem \ref{th:main1}).
First consider Theorems \ref{th:main1}.\ref{th:main1PN} and \ref{th:main1}.\ref{th:main1PandN}.
Let $Th(\mathbb{N})+(\Pi^{0}_{k}(\mbox{{\rm P}}),\mbox{{\rm P}}\cup\mbox{{\rm N}})\mbox{{\rm -ID}}\vdash R_{\varphi}(n)$
for a positive operator $\varphi(X,x)$,
or
$Th(\mathbb{N})+(\Pi^{0}_{k}(\mbox{{\rm P}}),\mbox{{\rm P}}\land\mbox{{\rm N}})\mbox{{\rm -ID}}({\sf Acc})\vdash R_{\varphi}(n)$
for an {\sf Acc}-operator $\varphi(X,x)$.
By Embedding Lemmas \ref{th:embedregthm1} and \ref{th:embedregthm2},
we have $\mathcal{H}_{2}\vdash_{2}^{a}I_{\varphi}^{<\Omega}(n)$
for $0<a<\Omega\cdot\omega_{1+k}$.
Then by Collapsing Lemma \ref{lem:collaps}
we obtain 
$\mathcal{H}_{\omega^{\Omega+a}+1}\vdash^{\psi(\omega^{\Omega+a})}_{1}I_{\varphi}^{<\Omega}(n)$, which in turn yields
$\mathcal{H}_{\omega^{\Omega+a}+1}\vdash^{\psi(\omega^{\Omega+a})}_{1}I_{\varphi}^{<\psi(\omega^{\Omega+a})}(n)$
by Bounding Lemma \ref{lem:bounding}.
We conclude 
$|n|_{\varphi}<\psi(\omega^{\Omega+a})<\psi(\omega^{\Omega\cdot \omega_{1+k}})=\psi(\Omega^{\omega_{1+k}})=\vartheta(\Omega\cdot\omega_{1+k})$.
 
Second consider Theorem \ref{th:main1}.\ref{th:main11}.
Let 
$Th(\mathbb{N})+\Pi^{0}_{1}(\mbox{{\rm P}})\mbox{{\rm -ID}}({\sf Acc})\vdash R_{\varphi}(n)$
for an {\sf Acc}-operator $\varphi(X,x)$.
By Embedding Lemma \ref{th:embedregthm3}
we have $\mathcal{H}_{2}\vdash_{3}^{a}I_{\varphi}^{<\Omega}(n)$
for $0<a<\Omega\cdot\omega$.
Lemma \ref{lem:CE} yields $\mathcal{H}_{2}\vdash_{2}^{\omega^{a}}I_{\varphi}^{<\Omega}(n)$.
Collapsing Lemma \ref{lem:collaps} together with Bounding Lemma \ref{lem:bounding}
yields 
$\mathcal{H}_{\omega^{\Omega+\omega^{a}}+1}\vdash^{\psi(\omega^{\Omega+\omega^{a}})}_{1}I_{\varphi}^{<\psi(\omega^{\Omega+\omega^{a}})}(n)$
We conclude $|n|_{\varphi}<\psi(\omega^{\Omega+\omega^{a}})<\psi(\omega^{\omega^{\Omega\cdot \omega}})=\psi(\Omega^{\Omega^{\omega}})=
\vartheta(\Omega^{\omega})$.

\appendix

\section{Resolution}\label{sec:appendix}

In {\bf Case 3} of the proof of Collapsing Lemma \ref{lem:collaps},
the following two are given.

\begin{equation}
\renewcommand{\theequation}{\ref{eq:case3.0}}
\mathcal{H}_{\gamma}\vdash^{a_{0}}_{2}\Gamma, \bigwedge_{i}(C_{i}\lor \bar{D}_{i})
\end{equation}
\addtocounter{equation}{-1}
and
\begin{equation}
\renewcommand{\theequation}{\ref{eq:case3.1}}
\mathcal{H}_{\gamma}\vdash^{a_{0}}_{2} \bigvee_{i}(\bar{C}_{i}\land D_{i}),\Gamma
\end{equation}
\addtocounter{equation}{-1}
Let $n$ be the number of conjunctions/disjunctions.
From (\ref{eq:case3.0}) and (\ref{eq:case3.1}) we obtain the following by inversion.
For each $i<n$
\beqn\label{eq:case3.0r}
\mathcal{H}_{\gamma}\vdash^{a_{0}}_{2}\Gamma, C_{i}, \bar{D}_{i}
\eeqn
and for each partition $I,J$ of the index set $\{0,\ldots,n-1\}$ ($I\cup J=\{0,\ldots,n-1\}$, $I\cap J=\emptyset$)
\beqn\label{eq:case3.1r}
\mathcal{H}_{\gamma}\vdash^{a_{0}}_{2} \{\bar{C}_{i}\}_{i\in I}, \{D_{j}\}_{j\in J},\Gamma
\eeqn

Let $b_{m}=\gamma+\omega^{\Omega+a_{0}}\cdot m$ and 
$\beta_{m}=\psi(b_{m})$ for positive integers $m$.
From (\ref{eq:case3.0r}) and (\ref{eq:case3.1r}) together with inversion, Bounding Lemma \ref{lem:bounding}
and the induction hypothesis,
we see the following.
\beqn\label{eq:case3.0rs}
\fal m,k \left[
m<k \Rarw
\mathcal{H}_{b_{k}+1}\vdash^{\bet_{k}}_{1}\Gamma, C_{i}^{(\bet_{k})}, \bar{D}_{i}^{(\bet_{m})}
\right]
\eeqn
and
\beqn\label{eq:case3.1rs}
\fal \vec{m},\vec{k} \left[
\max\vec{m}<\min\vec{k} \Rarw
\mathcal{H}_{b_{k}+1}\vdash^{\bet_{k}}_{1}\{\bar{C}_{i}^{(\bet_{m_{i}})}\}_{i\in I}, \{D_{j}^{(\bet_{k_{j}})}\}_{j\in J},\Gam
\right]
\eeqn
where $(I,J)$ is a partition,
$\vec{m}=(m_{i})_{i\in I}$ and $\vec{k}=(k_{j})_{j\in J}$ are sequences of positive integers.
$k=\max\{\max\vec{k},1+\max\vec{m}\}$, where $\max\vec{k}:=0$ and $\min\vec{k}:=\ome$
when $\vec{k}$ is the empty sequence.

We show that there exists a $k$ such that $\mathcal{H}_{b_{k}+1}\vdash^{\bet_{k}}_{1}\Gam$.
Now the above case $n=2$ indicates us that our task is to find a (ground) resolution refutation
and a `correct' assignment of integers to \textit{occurrences} of `literals' $C_{i},\bar{C}_{i},D_{i},\bar{D}_{i}$
in the refutation.

Since the set $\Gam$ of positive formulas plays no r\^{o}le, let us suppress it.
Second we omit the operator controlled part $\mathcal{H}_{\gam}\vdash^{a}_{d}$ since we can recover it.
Third let us denote $C^{(\bet_{m})}$ by $C^{(m)}$.
Thus (\ref{eq:case3.0r}) and (\ref{eq:case3.1r})
are written as follows.

\beqn
\renewcommand{\theequation}{\ref{eq:case3.0r}}
C_{i}, \bar{D}_{i}
\eeqn
\addtocounter{equation}{-1}

\beqn
\renewcommand{\theequation}{\ref{eq:case3.1r}}
 \{\bar{C}_{i}\}_{i\in I}, \{D_{j}\}_{j\in J}
\eeqn
\addtocounter{equation}{-1}

By the completeness of the ground resolution,
there exists a resolution refutation from `clauses' (\ref{eq:case3.0r}), (\ref{eq:case3.1r})
using only the ground resolution rule:
\[
\infer{\Gam,\Del}
{
\Gam,\bar{E}
&
E,\Del
}
\]
We need to find an assignment of positive integers to occurrences of literals
for which the following hold:
\benu
\item
If an integer $m$ is assigned to the left cut formula, i.e., the occurrence $\bar{E}$ in the ground resolution rule,
then the right cut formula $E$ receives the same number $m$.
$
\infer{\Gam,\Del}
{
\Gam,\bar{E}^{(m)}
&
E^{(m)},\Del
}
$

\item
If the two occurrences of the same literal are linked in the ground resolution rule,
then these two occurrences receive the same number.
This means if a literal $F$ occurs in the left upper sequent $\Gam$ and receives $m$, then the occurrence
of $F$ in the lower sequent $\Gam,\Del$ receives $m$.
The same for occurrences in $\Del$ and in $\Gam,\Del$.
$\infer{F^{(m)},\Gam,\Del}
{
F^{(m)},\Gam,\bar{E}
&
E,\Del,F^{(m)}
}
$

\item
The assignment to leaf clauses in (\ref{eq:case3.0r}) and (\ref{eq:case3.1r})
has to enjoy the conditions in (\ref{eq:case3.0rs}) and (\ref{eq:case3.1rs}).
This means if we attach numbers $k,m$ to a leaf clause $C_{i},\bar{D}_{i}$,
then $k>m$ has to hold:
\beqn
\renewcommand{\theequation}{\ref{eq:case3.0rs}}
C_{i}^{(k)}, \bar{D}_{i}^{(m)} \Rarw k>m
\eeqn
\addtocounter{equation}{-1}
Also if we attach numbers $\vec{m}$, $\vec{k}$ to a leaf clause
$\{\bar{C}_{i}\}_{i\in I}, \{D_{j}\}_{j\in J}$, then $\max\vec{m}<\min\vec{k}$ has to hold:

\beqn
\renewcommand{\theequation}{\ref{eq:case3.1rs}}
\{\bar{C}_{i}^{(m_{i})}\}_{i\in I}, \{D_{j}^{(k_{j})}\}_{j\in J}
\Rarw \max\vec{m}<\min\vec{k} 
\eeqn
\addtocounter{equation}{-1}

\eenu
When a ground resolution derivation together with an assignment of positive integers
to occurrences of literals enjoys the above three conditions,
the derivation together with attached integers is said to be \textit{decorated} derivation.

\bprp\label{prp:decoration}
For each $j<n$, there exist a $k_{j}$ and a decorated derivation $\pi_{j}(n)$ of $\bar{D}_{j}^{(k_{j})}$
from clauses (\ref{eq:case3.0r}) and (\ref{eq:case3.1r}).
\eprp
Assuming Proposition \ref{prp:decoration}, the following is a desired decorated refutation.
\[
\infer{\Box}
{
 \infer*[\pi_{0}(n)]{\bar{D}_{0}^{(k_{0})}}{}
 &
 \cdots
 &
 \infer*[\pi_{n-1}(n)]{\bar{D}_{n-1}^{(k_{n-1})}}{}
 &
  D_{0}^{(k_{0})},\ldots,D_{n-1}^{(k_{n-1})}
 }
 \]
where $D_{0}^{(k_{0})},\ldots,D_{n-1}^{(k_{n-1})}$ is one of clauses in (\ref{eq:case3.1r}), i.e., $I=\emptyset$.

We show Proposition \ref{prp:decoration} by induction on $n\geq 1$.
The case $n=1$ is clear: for any $1\leq k<m$ the following is a decorated derivation of $\bar{D}_{0}^{(k)}$.
\[
\pi_{0}(1)=
\infer{\bar{D}_{0}^{(k)}}
{
C_{0}^{(m)},\bar{D}_{0}^{(k)}
&
\bar{C}_{0}^{(m)}
}
\]
Assume that Proposition \ref{prp:decoration} holds for $n\geq 1$, and
let $\pi_{j}(n)$ be decorated derivation of $\bar{D}_{j}^{(k_{j})}$ for each $j<n$.

\bprp\label{prp:dec}
\benu
\item\label{prp:dec-}
For each $p<n$ and $m\geq 1$, there exist a $k_{p}^{\prime}\geq k_{p}$ and a decorated derivation 
$\pi_{p}(n)*\bar{C}_{n}^{(m)}$ of $\bar{D}_{p}^{(k^{\prime}_{p})},\bar{C}_{n}^{(m)}$.

\item\label{prp:dec+}
For each $q<n$, there exist a $k$ and a decorated derivation 
$\pi_{q}(n)* D_{n}^{(m)}$ of $\bar{D}_{q}^{(k_{q})},D_{n}^{(k)}$.
\eenu
\eprp
Assuming Proposition \ref{prp:dec} the following with $m>k$
 is a desired decorated derivation $\pi_{p}(n+1)$ of $\bar{D}_{p}^{(k_{p})}$
for $p<n$.
{\small
\[
\infer{\bar{D}_{p}^{(k^{\prime}_{p})}}
{
C_{n}^{(m)},\bar{D}_{n}^{(k)}
&
\hskip -10mm
\infer{D_{n}^{(m)}}
 {
  \infer*[\pi_{0}(n)* D_{n}^{(k)}]{\bar{D}_{0}^{(k_{0})},D_{n}^{(k)}}{}
  &
  \hskip -10mm
  \cdots
  &
   \infer*[\pi_{n-1}(n)* D_{n}^{(k)}]{\bar{D}_{n-1}^{(k_{n-1})},D_{n}^{(k)}}{}
  &
  \hskip -5mm
  D_{0}^{(k_{0})},\ldots,D_{n-1}^{(k_{n-1})},D_{n}^{(k)}
 }
&
\hskip -5mm
\infer*[\pi_{p}(n)*\bar{C}_{n}^{(m)}]{\bar{D}_{p}^{(k^{\prime}_{p})},\bar{C}_{n}^{(m)}}{}
}
\]
}
A decorated derivation $\pi_{n}(n+1)$ of $\bar{D}_{n}^{(k^{\prime}_{n})}$ is obtained from $\pi_{0}(n+1)$ of $\bar{D}_{0}^{(k^{\prime}_{0})}$
by interchanging indices $0$ and $n$ in literals $C_{0},D_{0},C_{n},D_{n}$.
\\

\noindent
{\bf Proof} of Proposition \ref{prp:dec}.

First we show Proposition \ref{prp:dec}.\ref{prp:dec-}.
Let $p<n$ and $m\geq 1$.
In the decorated derivation $\pi_{p}(n)$ of $\bar{D}_{p}^{(k_{p})}$,
append the negative decorated literal $\bar{C}_{n}^{(m)}$ to each leaf clause (\ref{eq:case3.1r}),
and then append $\bar{C}_{n}^{(m)}$ to each clause occurring below a leaf clause (\ref{eq:case3.1r}).
The result may violate the condition (\ref{eq:case3.1rs})
when $\max\vec{m}<\min\vec{k}\leq m$ in a leaf $\{\bar{C}_{i}^{(m_{i})}\}_{i\in I}, \{D_{j}^{(k_{j})}\}_{j\in J}$.
To avoid this, given an $m$, raise all of the assigned integers to literals in $\pi_{p}(n)$ by $1+m$.
Each $E_{i}^{(r)}$ is replaced by $E_{i}^{(1+m+r)}$.
This results in a decorated derivation of $\bar{D}_{p}^{(k^{\prime}_{p})}$ with $k^{\prime}_{p}=1+m+k_{p}$.

Next consider Proposition \ref{prp:dec}.\ref{prp:dec+}.
Let $q<n$.
In the decorated derivation $\pi_{q}(n)$ of $\bar{D}_{q}^{(k_{q})}$,
append the positive \textit{undecorated} literal $D_{n}$ to each leaf clause (\ref{eq:case3.1r}),
and then append $D_{n}$ to each clause occurring below a leaf clause (\ref{eq:case3.1r}).
We need to find an integer $k$ such that if we assign the number $k$ to all occurrences of $D_{n}$,
then the result is a decorated derivation $\pi_{q}(n)* D_{n}^{(k)}$ of $\bar{D}_{q}^{(k_{q})},D_{n}^{(k)}$.
For each leaf clause $\{\bar{C}_{i}^{(m_{i})}\}_{i\in I}, \{D_{j}^{(k_{j})}\}_{j\in J}$ in $\pi_{q}(n)$,
$\max\vec{m}<k$ is required for the condition (\ref{eq:case3.1rs}).
It suffices for $k$ to be larger than every assigned number $m_{i}$ to
the negative literal $\bar{C}_{i}$.
Thus 
$k=1+\max(\max\vec{m}:\vec{m})$ suffices, where $\vec{m}$ ranges over every decorated leaf clause 
$\{\bar{C}_{i}^{(m_{i})}\}_{i\in I}, \{D_{j}^{(k_{j})}\}_{j\in J}$ in $\pi_{q}(n)$.
\eprf

\end{document}